\documentclass[10pt, reqno]{amsart}
\usepackage{graphicx} 

\usepackage{extarrows}
\usepackage{stmaryrd}
\usepackage{MnSymbol}
\numberwithin{equation}{section}

\usepackage{amsxtra}
\usepackage{color}
\usepackage{hyperref}
\usepackage{tikz-cd}
\usepackage[all]{xy}
\usetikzlibrary{patterns}
\usetikzlibrary{backgrounds}

\usepackage{mathtools}
\newcommand{\mathdash}{\relbar\mkern-9mu\relbar}

\let\oldtocsection=\tocsection

\let\oldtocsubsection=\tocsubsection

\let\oldtocsubsubsection=\tocsubsubsection

\renewcommand{\tocsection}[2]{\hspace{0em}\oldtocsection{#1}{#2}}
\renewcommand{\tocsubsection}[2]{\hspace{2em}\oldtocsubsection{#1}{#2}}
\renewcommand{\tocsubsubsection}[2]{\hspace{4.5em}\oldtocsubsubsection{#1}{#2}}

\makeatletter
\def\subsection{\@startsection{subsection}{2}
 \z@{.5\linespacing\@plus.7\linespacing}{-.5em}
 {\normalfont\bfseries}}
\def\subsubsection{\@startsection{subsubsection}{3}
 \z@{.5\linespacing\@plus.7\linespacing}{-.5em}
 {\normalfont\bfseries}}
\makeatother

\theoremstyle{plain}

\renewcommand{\theequation}{\arabic{section}.\arabic{equation}}

\newtheorem{theorem}{Theorem}[section]
\newtheorem{lemma}[theorem]{Lemma}
\newtheorem{definition-theorem}[theorem]{Definition-Theorem}
\newtheorem{proposition}[theorem]{Proposition}
\newtheorem{corollary}[theorem]{Corollary}
\newtheorem{definition}[theorem]{Definition}
\newtheorem{example}[theorem]{Example}
\newtheorem{remark}[theorem]{Remark}
\newtheorem{notation}[theorem]{Notation}
\newtheorem{assumption}[theorem]{Assumption}
\newtheorem{lemma-definition}[theorem]{Lemma-Definition}
\newtheorem{lemma-notation}[theorem]{Lemma-Notation}
\newtheorem{question}[theorem]{Question}
\newtheorem{remark-definition}[theorem]{Remark-Definition}
\newtheorem{notation-remark}[theorem]{Notation-Remarks}
\newtheorem{conjecture}[theorem]{Conjecture}
\newcommand \bth[1] { \begin{theorem}\label{t#1} }
\newcommand \ble[1] { \begin{lemma}\label{l#1} }

\newcommand \bpr[1] { \begin{proposition}\label{p#1} }
\newcommand \bco[1] { \begin{corollary}\label{c#1} }
\newcommand \bconj[1] { \begin{conjecture}\label{co#1} }
\newcommand \bde[1] { \begin{definition}\label{d#1}\rm }
\newcommand \bex[1] { \begin{example}\label{e#1}\rm }
\newcommand \bre[1] { \begin{remark}\label{r#1}\rm }
\newcommand \bnota[1] {\begin{notation}\label{n#1}\rm }
\newcommand \bas[1] { \begin{assumption}\label{a#1}\rm }
\newcommand \bld[1] { \begin{lemma-definition}\label{ld#1} }

\newcommand \bqu[1] { \begin{question}\label{q#1}\rm }

\newcommand {\eth} { \end{theorem} }
\newcommand {\ele} { \end{lemma} }

\newcommand {\epr} { \end{proposition} }
\newcommand {\econj} { \end{conjecture} }
\newcommand {\eco} { \end{corollary} }
\newcommand {\ede} { \end{definition} }
\newcommand {\eex} { \end{example} }
\newcommand {\ere} { \end{remark} }
\newcommand {\enota} { \end{notation} }
\newcommand {\eas} {\end{assumption}}
\newcommand {\eld}{ \end{lemma-definition} }

\newcommand {\equ} {\end{question}}
\newcommand \thref[1]{Theorem \ref{t#1}}
\newcommand \leref[1]{Lemma \ref{l#1}}
\newcommand \prref[1]{Proposition \ref{p#1}}
\newcommand \coref[1]{Corollary \ref{c#1}}

\newcommand \deref[1]{Definition \ref{d#1}}
\newcommand \exref[1]{Example \ref{e#1}}
\newcommand \reref[1]{Remark \ref{r#1}}

\newcommand \ldref[1]{Lemma-Definition \ref{ld#1}}

\def \RR {{\mathbb R}}         
\def \CC {{\mathbb C}}
\def \ZZ {{\mathbb Z}}
\def \NN {{\mathbb N}}
\def \TT {{\mathbb T}}

\def \QQ {{\mathbb Q}}

\def \AA {{\mathbb A}}

\def \XX {{\mathbb X}}

\def \Lcc {{L^{c,c^{-1}}}}
\def \Acc {{A^{c,c^{-1}}}}
\def \Xcc {{X^{c,c^{-1}}}}

\def \calA {{\mathcal{A}}}
\def \calX {{\mathcal{X}}}
\def \calT {{\mathcal{T}}}

\def \calD {{\calD}}
\def \FF {\mathbb{F}}
\def \tkt{{\begin{xy}(0,1)*+{t}="A",(10,1)*+{t'}="B",\ar@{-}^k"A";"B" \end{xy}}}
\def \tskt{{\begin{xy}(0,1)*+{t}="A",(14,1)*+{t'}="B",\ar@{-}^{\sigma^{-1}(k)}"A";"B" \end{xy}}}
\def \tkv{{\begin{xy}(0,1)*+{t}="A",(10,1)*+{v}="B",\ar@{-}^k"A";"B" \end{xy}}}
\def \tskv{{\begin{xy}(0,1)*+{t}="A",(14,1)*+{v}="B",\ar@{-}^{\sigma(k)}"A";"B" \end{xy}}}

\def \vkv{{\begin{xy}(0,1)*+{v}="A",(10,1)*+{v'}="B",\ar@{-}^k"A";"B" \end{xy}}}

\def \fX {{\mathfrak{X}}}
\def \calC {{\mathcal C}}

\def \calP {{\mathcal P}}
\def \calQ {{\mathcal Q}}

\def \calD {{\mathcal D}}

\def \calU {{\mathcal U}}

\def \wh {\widehat}

\def \calApos {{\calA_{>0}}}
\def \calXpos {{\calX_{>0}}}
\def \Aprinpos {{\calA^{\rm prin}_{>0}}}
\def \Aprin {{\calA^{\rm prin}}}
\def \psipos {\psi_{>0}}

\def \DTA {{\rm DT}_\calA}
\def \DTX {{\rm DT}_\calX}

\def \DTApos {{\rm DT}_\calApos}
\def \DTXpos {{\rm DT}_\calXpos}
\def \dtA {{\rm dt}_{\calA}}
\def \dtX {{\rm dt}_{\calX}}

\def \dtAa {{\rm dt}_{\calA, {\bf a}}}
\def \dtXb {{\rm dt}_{\calX, {\bf b}}}

\def \g {\mathfrak{g}}
\def \t {\mathfrak{t}}

\def \ka {\kappa}
\def \hs {\hspace{.2in}}

\def \pp {\backslash \!\backslash}
\def \bSigma {\boldsymbol{\Sigma}}
\def \oc {\overline{c}}
\def \Accpos {A^{c, c^{-1}}_{>0}}
\def \Lccpos {L^{c, c^{-1}}_{>0}}
\def \Xccpos {X^{c, c^{-1}}_{>0}}

\title{Fixed points of DT transformations, cluster exponents and degrees of Weyl groups}
\author{Antoine de Saint Germain}
\address{
Department of Mathematics and New Cornerstone Science Laboratory   \\
The University of Hong Kong \\
Pokfulam Road               \\
Hong Kong}
\email{adsg96@hku.hk}
\author{Jiang-Hua Lu}
\address{
Department of Mathematics   \\
The University of Hong Kong \\
Pokfulam Road               \\
Hong Kong}
\email{jhlu@maths.hku.hk}

\begin{document}
\begin{abstract}
    Inspired by a recent work of Y. Mizuno, we show that the DT transformations on a cluster ensemble of finite type admit unique totally positive fixed points, and that the exponents of the linearizations of the  DT transformations at the fixed points
    are precisely the degrees of the Weyl group of the corresponding finite root system. 
\end{abstract}
\maketitle
\tableofcontents
\addtocontents{toc}{\protect\setcounter{tocdepth}{1}}

\section{Introduction and main results}
\subsection{Introduction} In the celebrated work \cite{FZ:II}, 
S. Fomin and A. Zelevinsky proved that cluster
algebras of finite type have the same Cartan-Killing classification as finite root systems. Direct connections between
the two theories are thus of great interest. For an indecomposable root system $\Phi$ of rank $n$, if
$\calA$ is a cluster $K_2$-variety (see $\S$\ref{ss:cluster-ensemble}) of the same Cartan-Killing type as $\Phi$, 
one set of invariants of $\Phi$ that appears\footnote{In \cite{FZ:associahedra}, Fomin and Zelevinsky refer 
to this appearance as a ``mystery".} in the combinatorics of $\calA$ are the degrees 
$d_1, d_2, \ldots, d_n$ of the Weyl group $W$ of $\Phi$, i.e., the (well-defined) degrees of algebraically independent
homogeneous generators of the ring of 
$W$-invariant polynomials on the vector space spanned by $\Phi$ 
(see \cite[Page 59]{HumphreysCoxGrp} and $\S$\ref{ss:cluster exponents} for detail). More precisely,  it is 
shown in \cite{FZ:associahedra} that 
\begin{equation}\label{eq:intro-0}
    \# \text{ of cluster variables of } \calA  = \sum_{i=1}^n d_i  \hs \mbox{and}  \hs
    \# \text{ of clusters of } \calA  = \prod_{i=1}^n \frac{d_i+{\bf h}}{d_i},
\end{equation}
where  ${\bf h}$ is the Coxeter number of $\Phi$. The identities in \eqref{eq:intro-0}  are
established in \cite{FZ:II, FZ:associahedra} through the connection between the exchange graph of $\calA$ 
and  a convex polytope $\Delta(\Phi)$, called a {\it generalized associahedron}, associated to the root system $\Phi$.
In particular, vertices of $\Delta(\Phi)$  correspond to clusters of $\calA$ and facets  to cluster variables
of $\calA$.

In this paper, we establish another connection between the invariants $d_1, d_2, \ldots, d_n$ of  $\Phi$
and the cluster variety $\calA$, this time through the {\it DT transformation} of $\calA$. 

\subsection{Statements of main results}
We will in fact work with 
cluster ensembles $(\calA,\calX,p)$ as defined by V.V. Fock and A.B. Goncharov \cite{FG:ensembles}, where $\calA$ is a cluster $K_2$-variety, $\calX$ is a cluster Poisson variety, and $p: \calA \to \calX$ is an {\it ensemble map} (see $\S$\ref{ss:cluster-ensemble}).
The {\it cluster Donaldson-Thomas transformations},
or {\it cluster DT transformations} for short, on a cluster ensemble  $(\calA,\calX,p)$
 are a pair of biregular automorphisms 
\[
    \DTA: \calA \longrightarrow \calA \quad \text{ and } \quad \DTX: \calX \longrightarrow \calX
\]
(see $\S$\ref{ss:DT-expon} for detail) which are a cluster-theoretic version of M. Kontsevich 
and Y. Soibelman's notion of DT transformations in the context of general 3-Calabi-Yau varieties
(see, for example, \cite{keller_2012, Nagao, Goncharov-Shen-DT}). The existence of cluster DT transformations has important consequences for the underlying cluster ensemble; for instance it implies that the Fock-Goncharov 
duality conjectures \cite{FG:ensembles} hold for the ensemble \cite{GHKK}. For this reason, cluster DT transformations have become an important object of study.

Let $\Phi$ be an indecomposable finite root system of rank $n$. Let $(\calA,\calX,p)$ be a
cluster ensemble of the same Cartan-Killing type as $\Phi$ ($\S$\ref{ss:finitetype}), and we will simply say that $(\calA,\calX,p)$ is of type $\Phi$. We also assume that both $\calA$ and $\calX$ have trivial coefficients. 
It is then well-known that  DT transformations for $(\calA,\calX,p)$ exist \cite{FZ:associahedra, GHKK} and are
periodic. More precisely, 
\begin{equation}\label{eq:DT-finite-order}
    (\DTA)^{{\bf h}+2} =  (\DTX)^{{\bf h}+2} = {\rm Id},
\end{equation}
where ${\bf h}$ is the Coxeter number of $\Phi$. In the literature, \eqref{eq:DT-finite-order} is commonly referred to as (a special case of) {\it Zamolodchikov periodicity}, postulated by A. Zamolodchikov in \cite{Zamolo:Y-system} and proved in full generality by B. Keller \cite{Keller:periodicity}.

Denote by $\calApos$ and $\calXpos$ the sets
of $\RR_{>0}$-valued points of $\calA$ and $\calX$ respectively. Then $\DTA$ preserves $\calApos$ and $\DTX$ preserves $\calXpos$ (see \S \ref{ss:DT-expon}). Let $\DTApos$ be the restriction of $\DTA$ to $\calApos $ and $\DTXpos$ the restriction of $\DTX$ to $\calXpos$. By \cite[Theorem B]{Ishibashi-fxpoint}, Zamolodchikov periodicity implies that
$\DTApos$ has fixed points in $\calApos$ and $\DTXpos$ has fixed points in $\calXpos$. We prove the uniqueness of 
such fixed points.

\medskip
{\it {\bf Theorem A} (\thref{fixedpoint}).
For any cluster ensemble $(\calA,\calX,p)$ of finite type, 
$\DTApos$ has a unique fixed point ${\bf a}$ in $\calApos$, and
$\DTXpos$ has a unique fixed point ${\bf b}$ in $\calXpos$. Moreover, $p({\bf a}) = {\bf b}$.
}
\medskip

As $\calApos$ and $\calXpos$ are parametrized by $(\RR_{>0})^n$ using any cluster, both $\calApos$ and $\calXpos$ are smooth real manifolds. With the fixed points ${\bf a} \in \calApos$ and ${\bf b} \in \calXpos$ as in Theorem A,
let $T_{{\bf a}} (\calApos)$ and $T_{{\bf b}} (\calXpos)$ be the respective
tangent space of $\calApos$ at ${\bf a}$ and of $\calXpos$ at ${\bf b}$, and set
\begin{equation}\label{eq:dt-A-X}
\dtA = {\rm d}_{\bf a}(\DTApos) \in {\rm End}_\RR (T_{{\bf a}} (\calApos)) \hs \mbox{and} \hs 
\dtX = {\rm d}_{\bf b}(\DTXpos) \in {\rm End}_\RR (T_{{\bf b}} (\calXpos)).
\end{equation}
The main results of our paper are as follows.

\medskip
\noindent
{\it {\bf Theorem B.} For any indecomposable finite root system $\Phi$ with Coxeter number ${\bf h}$ 
and any cluster ensemble $(\calA, \calX, p)$ 
of type $\Phi$, the following hold:

1) (\prref{:matrix-realization-linearDT}).  The two linear maps $\dtA$ and $\dtX$ have the same 
matrix representation in suitable bases of $T_{{\bf a}} (\calApos)$ and $T_{{\bf b}} (\calXpos)$.
In particular, $\dtA$ and $\dtX$ have the same characteristic polynomial, which we denote as $\calP_\Phi(x) = \det(x-\dtA)= \det(x-\dtX) \in \RR[x]$;

2) (\thref{:expnts-conjecture}). With $\zeta = e^{\frac{2\pi \sqrt{-1}}{{\bf h}+2}}$, one has 
\[
\calP_{\Phi} (x) = \prod_{j = 1}^n \left(x - \zeta^{d_j}\right),
\]
where $2 = d_1 \leq d_2 \leq \cdots \leq d_n = {\bf h}$ are the degrees of the Weyl group of $\Phi$.

3) (\prref{:degrees-and-heights}). The polynomial $\calP_{\Phi}(x)$  is the unique polynomial in $x$ satisfying
\begin{equation}\label{eq:char-poly-factor-intro}
    \left(\frac{x^{{\bf h}+2}-1}{x-1}\right)^n = \calP_\Phi (x)\prod_{\beta \in \Phi} (x - \zeta^{{\rm ht}(\beta)}),
\end{equation}
where ${\rm ht}: \Phi \rightarrow \ZZ$ is the height function associated to any choice of simple roots of $\Phi$. 
}
\medskip

When $\Phi$ is simply-laced, identity \eqref{eq:char-poly-factor-intro} implies a special case of
a conjecture of Mizuno \cite[Conjecture 3.8]{Mizuno:Y-sys-exponents} (see \reref{:relation-with-Y-sys-exponents} and \coref{:Mizuno-conj} for more detail), and Theorem B, 2) was recently used by A. King and T. Palaniappan \cite{King-2} to show that the multiplicity of $-1$ as an eigenvalue of $\dtA$ coincides with the dimension of the set of fixed points of $(\DTApos)^2$.

Let $G$ be a connected and simply connected complex Lie group with root system $\Phi$, and let
$c$ be a Coxeter element in the Weyl group $W$. A
geometric model of the cluster ensemble $(\calA, \calX, p)$ is provided by the reduced double Bruhat cell
$\Lcc \subset G$, and the  DT transformations of $(\calA, \calX, p)$  can be realized 
\cite[Theorem 1.1]{Daping-DT-doubleBruhat}
via the Fomin-Zelevinsky twist $\psi: \Lcc \to \Lcc$ introduced in \cite{FZ:double}. Let $L^{c, c^{-1}}_{>0}\subset \Lcc$ be
the totally positive part
of $\Lcc$ defined by Lusztig \cite{Lusztig1994}.
In $\S$\ref{s:Lie} we describe the fixed point
${\bf a} \in \calApos \subset L^{c, c^{-1}}_{>0}$, and we show that ${\bf a}$ is also the unique
fixed point of $\psi$ in $L^{c, c^{-1}}_{>0}$. We further prove in \thref{:Lie} that the 
characteristic polynomial of the linear map 
\[
d_{\bf a} \left(\psi|_{L^{c, c^{-1}}_{>0}}\right):\;\; T_{\bf a} \left(L^{c, c^{-1}}_{>0}\right) \longrightarrow 
T_{\bf a} \left(L^{c, c^{-1}}_{>0}\right)
\]
is the product polynomial $\calP_{\Phi}(x)\calP_c(x)$, where $\calP_c(x)$ is the characteristic polynomial of $c \in W$ as a linear operator on the Lie algebra of a maximal torus of $G$.


\subsection{Cluster exponents} Motivated by Theorem B, we introduce the notion of 
 {\it cluster exponents} for general cluster ensembles.
We first recall a definition from linear algebra.

\bde{:linear-expo-intro}
{\rm 
For an $n$-dimensional real vector space $E$ and $T \in {\rm GL}(E)$ such that
$T^{\bf q} = {\rm Id}_E$
for some integer ${\bf q} \geq 1$, the {\it exponents of $T$ with respect to ${\bf q}$}
are the integers $0 \leq m_1 \leq \cdots \leq m_n\leq {\bf q}-1$ such that 
\[
\det(x-T) = \prod_{i=1}^n \left(x - e^{\frac{2\pi\sqrt{-1}}{{\bf q}}m_i}\right).
\]
\hfill $\diamond$
}
\ede

\bex{ex:c}
{\rm 
When $V$ is the vector space spanned by a finite root system $\Phi$, the order of any Coxeter element
$c$ in the Weyl group $W$ of $\Phi$ is the Coxeter number ${\bf h}$, 
and it is 
well-known  \cite[Theorem 3.19]{HumphreysCoxGrp} that 
the exponents of  $c \in W$ with respect to ${\bf h}$, as a linear operator on $V$, 
are $d_1-1,  d_2-1, \cdots, d_n-1$, where again $2 = d_1 \leq d_2 \leq  \cdots 
\leq d_n = {\bf h}$
are the degrees of $W$. 
}
\eex

Suppose now that a cluster ensemble $(\calA, \calX, p)$, {\it not necessarily of finite type}, has
trivial coefficients and admits
DT transformations $\DTA: \calA \to \calA$ and $\DTX: \calX \to \calX$ (see $\S$\ref{ss:DT-expon} for detail).
Suppose that ${\bf a} \in \calApos$ is a fixed point of $\DTApos: \calApos \to \calApos$. Then 
${\bf b} :=p({\bf a}) \in \calXpos$ is a fixed point of $\DTXpos: \calXpos \to \calXpos$. 
Define their respective linearizations 
\[
\dtAa = {\rm d}_{\bf a}(\DTApos) \in {\rm End}_\RR (T_{{\bf a}} (\calApos)) \hs \mbox{and} \hs 
\dtXb = {\rm d}_{\bf b}(\DTXpos) \in {\rm End}_\RR (T_{{\bf b}} (\calXpos)).
\]
The following Theorem C is part of \thref{:same-char}.

\medskip
\noindent
{\bf Theorem C.}  {\it The linear maps $\dtAa$ and $\dtXb$ have the same characteristic polynomial.}

\medskip
Theorem C is proved by lifting the DT transformations on $\calA$ and $\calX$ to a 
{\it twisted automorphism $\psi$ 
of DT-type} \cite[\S 5.1]{kimura-qin-wei} on an $\AA$-cluster variety $\Aprin$
with principal coefficients called a
{\it principal extension} of $\calA$ (see \deref{:prin-extension}),
and by computing the linearization of $\psi$ at ${\bf a}$ in two coordinate systems. For any 
such a lifting $\psi:\Aprin\to\Aprin$, let $E$ be the unique $n \times n$ integer matrix such that 
$\psi^*({\bf p}) = {\bf p}^E$, where ${\bf p}=(p_1, \ldots, p_n)$ are the frozen variables of $\Aprin$.
By viewing ${\bf a} \in \calApos$ as a fixed point of $\psi$ in the positive part
$\calA^{\rm prin}_{>0}$ of $\Aprin$, we prove in 
\thref{:same-char} that  the characteristic polynomial of the linear map
\[
d_{\bf a}\left(\psi|_{\calA^{\rm prin}_{>0}}\right): \;\; T_{\bf a}(\calA^{\rm prin}_{>0})
\longrightarrow T_{\bf a}\left(\calA^{\rm prin}_{>0}\right)
\]
is equal to $\calP(x) \calP_E(x)$, where
$\calP(x)$ is the common characteristic polynomial of $\dtAa$ and $\dtXb$, and
$\calP_E(x) = \det(x-E)$ is the characteristic polynomial of $E$. As a special case 
we obtain, see \thref{:Lie} for detail, the statement on the characteristic polynomial of
the linearization of the 
Fomin-Zelevinsky twist
on the reduced double Bruhat cell $\Lcc$.

\medskip
If a cluster ensemble $(\calA, \calX, p)$ admits DT transformations 
$\DTA: \calA \to \calA$ and $\DTX: \calX \to \calX$, 
 the {\it synchronicity phenomenon} in cluster theory \cite{Nakanishi-Synchro} implies that 
if $\DTA$ is periodic with period ${\bf q}$, so is $\DTX$.
Motivated by the finite type case, we introduce the following definition.

\bde{:cluster-expo-intro}
{\rm Suppose that $(\calA, \calX, p)$ is a cluster ensemble (with trivial coefficients)  admitting 
DT transformations $\DTA: \calA \to \calA$ and $\DTX: \calX \to \calX$, and suppose that
$\DTA$ and $\DTX$ are periodic with period ${\bf q}$.
Suppose also that ${\bf a} \in \calApos$ is a unique fixed point of 
$\DTApos: \calApos \to \calApos$. We call the 
exponents of $\dtAa \in {\rm End}_{\RR}(T_{\bf a}(\calApos))$  with respect to ${\bf q}$  the {\it cluster exponents} of $(\calA, \calX, p)$
with respect to ${\bf q}$.
\hfill $\diamond$ 
}
\ede

\bre{:expon}
{\rm 
The assumption that ${\bf a} \in \calApos$ is a unique fixed point of $\DTApos$ in
$\calApos$ implies that $\dtA$ does not have 
$1$ as an eigenvalue (see \reref{:fpIsLefschetz}). The cluster exponents $m_1 \leq \cdots \leq m_n$ 
of $(\calA, \calX, p)$ with respect to ${\bf q}$ 
then satisfy  
$1\leq m_1 \leq \cdots \leq m_n \leq {\bf q}-1$. Since the characteristic polynomial of $\dtA$ has
real coefficients, one has $m_i + m_{n+1-i} = {\bf q}$ for every $i \in [1, n]$. 
It follows that 
\begin{equation}\label{eq:sum-mi}
\sum_{i=1}^n m_i= \frac{n{\bf q}}{2}.
\end{equation}
\hfill $\diamond$
}
\ere

For a cluster ensemble $(\calA, \calX, p)$ of type $\Phi$, where $\Phi$ is 
a finite dimensional root system 
with Coxeter number ${\bf h}$,  
Theorem B says that the 
cluster exponents of $(\calA, \calX, p)$ with respect to ${\bf h} + 2$ coincide with the degrees 
$2 =d_1 \leq \cdots \leq d_n = {\bf h}$ of the Weyl group of $\Phi$. 
Equation \eqref{eq:sum-mi} in this case gives\footnote{Once we know that $1 = d_2-1\leq \cdots \leq d_n-1 = {\bf h}-1$
are the exponents of a Coxeter element $c \in W$ with respect to ${\bf h}$, identity 
\eqref{eq:sum-di} also follows by applying \eqref{eq:sum-mi} to $c$.}
\begin{equation}\label{eq:sum-di}
\sum_{i=1}^n d_i= \frac{n({\bf h}+2)}{2} = \frac{n {\bf h}}{2} + n 
\end{equation}
where $\frac{n {\bf h}}{2} + n$, being the number of almost positive roots in $\Phi$,
is well-known  \cite{FZ:II} to be 
equal to the number of cluster
variables of $\calA$. What is less known, it seems to the authors, is that the same number 
$\frac{n {\bf h}}{2} + n$ is also equal to the number of {\it global $y$-variables
of $\calX$}, i.e., those $y$-variables of $\calX$ that are also regular functions on $\calX$ 
(see $\S$\ref{ss:cluster-ensemble}), a fact that is proved \cite[Lemma 7.1]{CdsgL:ringel-conj}, where 
an explicit bijection between the cluster variables of $\calA$ and the global $y$-variables of $\calX$ 
is constructed. By extending
\eqref{eq:intro-0}  to
\begin{align}\label{eq:intro-22}
   & \# \text{ of cluster variables of } \calA  = \# \text{ of global $y$-variables of } \calX = \sum_{i=1}^n d_i,\\
    \label{eq:intro-11}
   & \# \text{ of clusters of } \calA  = \# \text{ of clusters of } \calX = \prod_{i=1}^n \frac{d_i+{\bf h}}{d_i},
\end{align}
and by interpreting the $d_i$'s appearing in \eqref{eq:intro-22} and \eqref{eq:intro-11} as cluster exponents
of the ensemble $(\calA, \calX, p)$, 
we see that 
cluster exponents encode important information about cluster ensembles.

For other examples of cluster ensembles  admitting periodic DT transformations, possibly with frozen variables, 
such as those associated to
the Grassmannian \cite{Weng:Grass}, it would be very interesting to look at the fixed points of the DT transformations and formulate and 
compute their cluster exponents. More generally, one may define and study exponents 
associated to periodic elements of the cluster modular group (see \cite{Ishibashi-fxpoint}, in particular \cite[Theorem 2.2]{Ishibashi-fxpoint}), exponents associated to periodic automorphisms (see \cite{Karp:cycle-map} and \cite[Appendix A]{GL:twist-richardson}) or exponents associated to ``mutation loops". For examples of the latter, see Y. Mizuno's work \cite{Mizuno:Y-sys-exponents} on {\it exponents of $Y$-systems}, 
which inspired our notion of cluster exponents. 

\subsection{Main steps in the proofs of Theorem A and Theorem B}
Let $\Phi$ be an indecomposable finite root system with Coxeter number ${\bf h}$, and let
$A = (a_{i, j})_{i, j \in [1, n]}$ be the Cartan matrix associated to a choice of simple roots in $\Phi$. 
In $\S$\ref{ss:finitetype}, we first recall the construction of a cluster ensemble 
$(\calA, \calX, p)$ of type $\Phi$, and we recall the explicit 
description of the DT transformations $\DTA: \calA \to \calA$ and $\DTX: \calX \to \calX$.
Using a certain cluster coordinate chart of $\calA$ to identify $\calApos$ with 
$(\RR_{>0})^n$, we prove in 
\prref{:ab-equations} that 
$\DTApos$-fixed points in $\calApos$ are in bijection with solutions  $(a_1, \ldots, a_n) \in (\RR_{>0})^n$ to the
equations 
\begin{equation}\label{eq:a-intro}
        a_i^2 = 1 + \prod_{j\neq i}a_j^{-a_{j,i}}, \quad i \in [1,n].
\end{equation}
Similarly, using certain (not necessarily cluster) coordinate chart of $\calX$ consisting of global
$y$-variables of $\calX$, we show in \prref{:ab-equations}  that 
$\DTXpos$-fixed points in $\calXpos$ are in bijection with solutions  $(b_1, \ldots, b_n) \in (\RR_{>0})^n$ to the
equations 
\begin{equation}\label{eq:b-intro}
 b_i^2 = \prod_{j\neq i}(1 + b_j)^{-a_{j,i}}, \quad i \in [1,n].
\end{equation}
Moreover, if $(a_, \ldots, a_n) \in (\RR_{>0})^n$ is a solution to \eqref{eq:a-intro}, then $(b_1, \ldots, b_n)
=(a_1^2-1, \ldots, a_n^2-1)$
is a solution to \eqref{eq:b-intro}. 
That the systems \eqref{eq:a-intro} and \eqref{eq:b-intro} admit unique positive solutions is known in the literature
 on $T$-systems, $Y$-systems, and $Q$-systems (see the proof of \leref{:fixpt-lemma2}). 
 Thus Theorem A holds. We have solved the equations \eqref{eq:a-intro} explicitly 
case-by-case for every indecomposable Cartan matrix of finite type, and the answers are given in 
Table \ref{table:A} in Appendix A. 

To prove Theorem B, we introduce, for 
any $\lambda = (\lambda_1,\ldots, \lambda_n) \in (\RR^\times)^n$,  the matrices 
\begin{equation}\label{eq:LU-intro}
U_{\lambda} = \left(\begin{array}{ccccc} \lambda_1 & a_{1,2} & a_{1, 3} & \cdots & a_{1, n}\\
0 & \lambda_2 & a_{2,3} & \cdots & a_{2, n}\\
0 & 0 & \lambda_3 & \cdots & a_{3, n}\\
\cdots & \cdots & \cdots & \cdots & \cdots\\
0 & 0 & 0 & \cdots & \lambda_n\end{array}\right) \hs \mbox{and} \hs
L_{\lambda} = 
\left(\begin{array}{ccccc} \lambda_1 & 0 & 0 & \cdots & 0\\
a_{2, 1} & \lambda_2 & 0 &\cdots & 0\\
a_{3,1} & a_{3,2} & \lambda_3 & \cdots & 0\\
\cdots & \cdots & \cdots & \cdots & \cdots\\
a_{n, 1} & a_{n, 2} & a_{n, 3} & \cdots & \lambda_n
\end{array}\right).
\end{equation}
With the unique positive solutions $(a_1, \ldots, a_n)$ to \eqref{eq:a-intro} and $(b_1, \ldots, b_n)$ to \eqref{eq:b-intro},
set 
\begin{equation}\label{eq:ci-intro}
\kappa = (\kappa_1,\ldots, \kappa_n) \in (\RR_{>1})^n\hs \mbox{with} \hs \kappa_i = \frac{a_i^2}{a_i^2-1} = \frac{b_i+1}{b_i} = \prod_{j=1}^n a_j^{a_{j, i}}, \hs i \in [1, n].
\end{equation}
We prove in \prref{:matrix-realization-linearDT} that 
the two linear maps $\dtA$ and $\dtX$ 
are represented, in suitable bases of $T_{{\bf a}} (\calApos)$ and $T_{{\bf b}} (\calXpos)$ respectively,  
by the same matrix
\begin{equation}\label{eq:Mc}
    M_{\kappa}: = -L_{\kappa} U_{\kappa}^{-1}.
\end{equation}
Denote the common  characteristic polynomial of $\dtA$ and $\dtX$ by $\calP_\Phi(x)$. Then  
\begin{equation}\label{eq:calP-intro}
\calP_\Phi(x) = \frac{1}{\kappa_1\kappa_2\cdots \kappa_n} \det (xU_{\kappa}+L_{\kappa}).
\end{equation}

\bre{:twist-and-cox-elmt} The matrix $M_{\bf 1}$, where ${\bf 1} = (1,\ldots , 1)$, 
represents a Coxeter element in the basis of fundamental weights 
(see \eqref{eq:M-one}). Thus $\calP_\Phi(x)$ can be regarded a deformation of 
the characteristic polynomial of any Coxeter element as a linear operator on 
the vector space spanned by $\Phi$. 
\hfill $\diamond$
\ere

Let again  $2 = d_1 \leq \cdots \leq d_n= {\bf h}$ be the degrees of the Weyl group of $\Phi$, and
for $\zeta = e^{\frac{2\pi \sqrt{-1}}{{\bf h}+2}}$, define
\begin{equation}\label{eq:calD-intro}
\calD_\Phi(x) = \prod_{j = 1}^n \left(x - \zeta^{d_j}\right).
\end{equation}
Using a theorem of Kostant relating 
$d_1, \ldots, d_n$ and the heights of roots with 
respect to simple roots in $\Phi$, we prove in \prref{:degrees-and-heights} that the polynomial $\calD_\Phi(x)$ satisfies
\eqref{eq:char-poly-factor-intro} with $\calD_\Phi(x)$ replacing $\calP_\Phi(x)$. 

In Appendix A, we use the explicit values of the Coxeter number ${\bf h}$ and 
of  $d_1, \ldots, d_n$ given in 
\cite[p. 59]{HumphreysCoxGrp} (see Table \ref{table:d}) to compute the polynomial 
$\calD_\Phi(x)$ for each $\Phi$ and state the results in \leref{:D-poly}. 
In Appendix B, 
we use \eqref{eq:calP-intro} to explicitly compute the polynomials $\calP_\Phi(x)$ 
and show directly that $\calP_\Phi(x) = \calD_\Phi(x)$ for each $\Phi$, thereby proving 2) and 3) of Theorem B.

\bre{:open}
{\rm
While it would be very interesting to find a uniform proof that $\calP_\Phi(x)=\calD_\Phi(x)$ for every 
indecomposable finite root system $\Phi$, it is also desirable 
to separately compute these polynomials explicitly and compare them, as we do in Appendices A and B. 
When  $\Phi$ can be obtained by $\Phi'$ via folding, our calculations 
also reveal, without appealing to the technique of folding in cluster algebras, that 
$\calP_{\Phi}$ is a factor of $\calP_{\Phi'}$. See Appendix B for details. 
\hfill $\diamond$
}
\ere



\subsection{Acknowledgments} 
The authors would like to thank Y. Mizuno 
for very useful comments and Zihang Liu and Peigen Cao for discussions. Zihang Liu in particular helped us in proving
 2) of Theorem B for 
type $C_n$ using a method different from the one presented in this paper. The second author would also like to thank Tuen Wai Ng, Matthew Pressland and Tsukasa Ishibashi for stimulating discussions on cluster exponents. A. de Saint Germain's research  has been partially supported by the New Cornerstone Science Foundation through the New Cornerstone Investigator Program awarded to Professor Xuhua He.
J.-H. Lu's research has been partially supported by the Research Grants Council 
(RGC) of the Hong Kong SAR, China (GRF 17306621 and 17306724). 

\subsection{Notational conventions}\label{ss:nota} The set of all positive real numbers is denoted by $\RR_{>0}$.
The transpose of a matrix $M$ is denoted as
$M^T$.
For an integer $n\geq 1$, we write $[1,n] = \{1,\ldots, n\}$. Elements in $\ZZ^n$ are regarded as 
column vectors. For $v \in \ZZ^n$, we write $v \geq 0$ if all its entries are non-negative, and 
similarly for $v \leq 0$. 
The standard basis of $\ZZ^n$ is denoted as $\{e_1, \ldots, e_n\}$. 

For any $\sigma$ in the permutation group $S_n$ on the set $[1, n]$, 
let $\overline{\sigma}\in {\rm GL}(n, \ZZ)$ be the 
$n \times n$  matrix 
\begin{equation}\label{eq:sigma-perm}
\overline{\sigma} = (e_{\sigma(1)}, \; e_{\sigma(2)}, \; \ldots, \; e_{\sigma(n}),
\end{equation}
so the map $S_n \to {\rm GL}(n, \ZZ), \sigma \mapsto \overline{\sigma}$, is a group homomorphism.
 If $V = (v_{i, j})$ is an $m \times n$ matrix and if $\tau \in S_m$
and $\sigma \in S_n$, writing $\overline{\tau}^{\, -1} \,V \,\overline{\sigma} = (v_{i, j}^\prime)$, one then has
$v_{i, j}^\prime = v_{\tau(i), \sigma(j)}$ for $i \in [1, m], j \in [1, n]$.

Let $A$ be a commutative ring with  group of units $A^\times$, and let
$a = (a_1, \ldots, a_n) \in (A^\times)^n$. We write $a^{-1} = (a_1^{-1}, \ldots, a_n^{-1})$.
For $V = (v_1, \ldots, v_n)^T \in \ZZ^n$ we write 
$a^V = a_1^{v_1} \cdots a_n^{v_n} \in A^\times$, and if 
$V$ is an $n \times m$ integral matrix with columns $V_1, \ldots, V_m$, we write
\begin{equation}\label{eq:aV}
a^V = (a^{V_1}, \ldots, a^{V_m}) \in (A^\times)^m.
\end{equation}
For $a  \in (A^\times)^n$ and two integer matrices $V$ and $W$ of sizes $n \times m$ and  $m \times k$
respectively, one then has
\begin{equation}\label{eq:a-VW}
a^{VW} = (a^V)^W \in (A^\times)^k.
\end{equation}
As a special case of \eqref{eq:aV}, for $\sigma \in S_n$ and  $a = (a_1, \ldots, a_n) \in (A^\times)^n$,
one has
\begin{equation}\label{eq:a-sigma}
a^{\overline{\sigma}} = a \, \overline{\sigma} = (a_{\sigma(1)}, \, \ldots, \, a_{\sigma(n)}).
\end{equation}
If $M$ is a left $A$-module, for $a = (a_1, \ldots, a_n) \in A^n$ and $\xi = (\xi_1, \ldots, \xi_n)\in M^n$, 
we write
\[
a\xi = (a_1\xi_1, \ldots, a_n\xi_n) \in M^n.
\]
In particular, for 
$a, b \in (A^\times)^n$, for any $n \times m$ integer matrix $V$ and for any $\sigma \in S_n$, we have
\begin{equation}\label{eq:ab-V-sigma}
(ab)^V = a^V b^V = (A^\times)^m, \hs (ab)^{\overline{\sigma}} = (ab)\overline{\sigma} = (a\overline{\sigma})
(b \overline{\sigma}).
\end{equation}

As an example, consider $(\RR_{>0})^n$ as a smooth manifold with coordinate functions
$x = (x_1, \ldots, x_n)$. Then for $V \in \ZZ^n$, the differential of the function
$x^V$ on $(\RR_{>0})^n$ satisfies
\begin{equation}\label{eq:d-xV-1}
x^{-V}d(x^V) = \left(\frac{dx_1}{x_1}, \, \ldots, \, \frac{dx_n}{x_n}\right) V = (x^{-1}dx)V.
\end{equation}
Consequently, if $V$ is an $n \times n$ integer matrix and $\phi: (\RR_{>0})^n \to (\RR_{>0})^n$ is given by
$\phi(x) = x^V$, then 
\begin{equation}\label{eq:d-xV-2}
\phi^*\left(\frac{dx_1}{x_1}, \, \ldots, \, \frac{dx_n}{x_n}\right) =
\left(\frac{dx_1}{x_n}, \, \ldots, \, \frac{dx_1}{x_n}\right) V, \hs \mbox{or}\hs
\phi^*(x^{-1}dx) = (x^{-1}dx) V.
\end{equation}

\section{DT transformations and cluster exponents}\label{s:DT-expon}
\subsection{Cluster ensembles}\label{ss:cluster-ensemble}
We recall some background on cluster ensembles and set up notation. 
All the material here is standard, and we refer to 
\cite{fz_2002,FZ:ClusterIV, FG:ensembles, GHK, GHKK} for more details.

For an integer $n \geq 1$, an $n \times n$ integer matrix $B$ is called a
{\it mutation matrix} if it is 
skew-symmetrizable,   i.e., if $DB$ is skew-symmetric for some diagonal matrix $D$ with positive 
integer entries on the diagonal. For another integer $m \geq 0$, an 
{\it extended mutation matrix} of size $(n+m) \times n$ is any 
$(n+m) \times n$ integer matrix $\widetilde{B}$, written as
$\widetilde{B}=(B \backslash \!\backslash  P)$, 
such that $B$ is skew-symmetrizable, where $B$ and $P$ 
are the sub-matrices of $\widetilde{B}$ respectively formed by the first $n$ rows and the 
last $m$ rows. The mutation matrix $B$ is called the {\it principal part} of the extended mutation matrix 
$\widetilde{B}$.

Let $n \geq 1$ and $m \geq 0$ be integers and let $\FF_{n+m}$ be a field 
extension of $\QQ$ of  pure transcendental degree $n+m$.
A {\it labeled seed in $\FF_{n+m}$} of rank $n$ is
a pair 
\begin{equation}\label{eq:Sigma}
\Sigma= ((x_1, \ldots, x_n, p_1, \ldots, p_m), \, \widetilde{B} = (b_{i,j})),
\end{equation}
where $(x_1, \ldots, x_n, p_1, \ldots, p_m)$ is an ordered set of free 
transcendental generators of $\FF_{n+m}$ over $\QQ$, and $\widetilde{B}$ is an extended
mutation matrix of size 
$(n+m) \times n$. Given $k \in [1, n]$, the {\it $\AA$-seed mutation} in direction $k$ of $\Sigma$ as in \eqref{eq:Sigma}  is 
the labeled seed
$\mu_k^\AA(\Sigma)=((x_1^\prime, \ldots, x_n^\prime, p_1, \ldots, p_m), \, \widetilde{B}^\prime = (b_{i, j}^\prime))$, where  
\begin{align}\label{eq:X-mut}
 &x_i' = \begin{cases}
  x_k^{-1}(\prod_{j=1}^{n+m} x_j^{[b_{j,k}]_+} +\prod_{j=1}^{n+m} x_j^{[-b_{j,k}]_+}) ,& i = k; \\
  x_i,& \text{otherwise};
  \end{cases}\\
\label{eq:matrix-mut}
		&b_{i,j}' = \begin{cases}
			-b_{i,j}, & \text{ if } i = k \text{ or } j = k, \\
			b_{i,j} + [b_{i,k}]_+[b_{k,j}]_+ - [-b_{i,k}]_+[-b_{k,j}]_+, &\text{ otherwise},
		\end{cases}
\end{align}
where $[x]_+ = \max (0,x)$ for $x \in \ZZ$ and $x_{n+j}= p_j$ for $j \in [1,m]$. One has 
$(\mu_k^\AA)^2(\Sigma)= \Sigma$ for every $k \in [1, n]$. 

Let $\TT_n$ be an 
$n$-regular tree with edges labeled by  the set $[1, n]$. 
A vertex of $\TT_n$ is denoted by $t$, and we write $t \in \TT_n$.
 If $t, t' \in \TT_n$ are
joined by an edge labeled by $k \in [1, n]$, we write \tkt. 

An {\it $\AA$-cluster pattern} in $\FF_{n+m}$ of rank $n$ is an assignment 
\begin{equation}\label{eq:bold-Sigma-A}
\boldsymbol{\Sigma}^\AA = \{\Sigma_t^\AA =({\bf x}_t=(x_{t; 1}, \ldots, x_{t; n}), {\bf p}=(p_1, \ldots, p_m)), \widetilde{B}_t = (B_t \pp P_t)) \}_{t \in \TT_n}
\end{equation}
of a labeled seed $\Sigma_t^\AA$ of rank $n$ in $\FF_{n+m}$ to each $t \in \TT_n$ such that 
$\mu_k^\AA(\Sigma_t^\AA) = \Sigma_{t'}^\AA$ whenever \tkt.
Given  $\boldsymbol{\Sigma}^\AA$ as in \eqref{eq:bold-Sigma-A},  
for each $t \in \TT_n$ one has the $\QQ$-rational torus 
$\calT_t={\rm Spec}(\QQ[{\bf x}_t^{\pm 1}, {\bf p}^{\pm 1}])$,
and for any $t, t' \in \TT_n$ one has a unique  birational morphism $\phi_{t', t}^\AA$
from $\calT_t$ to $\calT_{t'}$ which is the composition of the one-step mutations in \eqref{eq:X-mut}
along the unique path in $\TT_n$ from $t$ to $t'$. The 
{\it $\AA$-cluster variety}, or the {\it cluster-$K_2$ variety}, associated to $\boldsymbol{\Sigma}^\AA$, is the 
scheme $\calA$ over $\QQ$ obtained by gluing the 
tori $\{\calT_t\}_{t \in \TT_n}$ via the birational morphisms $\{\phi_{t', t}^\AA\}_{t, t' \in \TT_n}$.
The algebra of regular functions on $\calA$ is then the upper cluster algebra
(see \cite{GHK} for detail)
\[
\calU(\calA) = \bigcap_{t \in \TT_n} \QQ[{\bf x}_t^{\pm 1}, {\bf p}^{\pm 1}].
\]
Elements in $\mathfrak{X}(\calA):=\{x_{t; i}: t \in \TT_n, \, i \in [1, n]\}$ are called
{\it cluster variables} of $\calA$, and $p_1, \ldots, p_n$
are called the {\it frozen variables} of $\calA$. 
Each ${\bf x}_t$, for $t \in \TT_n$,
is called a {\it cluster} of $\calA$, while  $({\bf x}_t, {\bf p})$ is called an {\it extended cluster} of $\calA$.
We denote the set of all clusters of $\calA$ by $\calC(\calA)$. The assignment $\{\widetilde{B}_t\}_{t \in \TT_n}$ is called the {\it matrix pattern} of $\calA$.
When $m = 0$, the cluster $\AA$-variety $\calA$ is said to have rank $n$ with {\it trivial coefficients}.

For $\XX$-cluster varieties, also called {\it cluster-Poisson varieties}, we only consider the case of trivial coefficients. Let thus $\FF_n$ be a field extension of $\QQ$ of pure transcendental degree $n$. 
For $k \in [1, n]$, 
the {\it $\XX$-seed mutation} in direction $k$ of a rank $n$ 
labeled seed 
$({\bf y}, B) = ((y_1, \ldots, y_n), B)$ in $\FF_n$  is the labeled seed 
$\mu_k^\XX({\bf y}, B) =((y_1^\prime, \ldots, y_n^\prime), B')$ in $\FF_n$, where 
\begin{equation}\label{eq:Y-mut}
y_{i}^\prime = \begin{cases}
			y_{k}^{-1}, &  i = k; \\
			y_{i} y_{k}^{[b_{k,i}]_+} (1+y_{k})^{-b_{k,i}}, &\text{otherwise},
		\end{cases}
\end{equation}
and $B' = (b_{i, j}^\prime)$ is obtained from $B = (b_{i, j})$ by \eqref{eq:matrix-mut} (for $m = 0$). 
An {\it $\XX$-seed pattern} in $\FF_n$ is an assignment 
\begin{equation}\label{eq:bold-Sigma-X}
\bSigma^\XX = \{\Sigma_t^\XX =({\bf y}_t = (y_{t; 1}, \ldots, y_{t; n}), \;B_t)\}_{t \in \TT_n}
\end{equation}
of a seed in $\FF_n$ of rank $n$ to each $t \in \TT_n$ such that
$\Sigma_{t'}^\XX = \mu_k^\XX(\Sigma_t^\XX)$ whenever \tkt. Similar to the case of 
$\AA$-cluster varieties, 
given $\bSigma^\XX$ as in \eqref{eq:bold-Sigma-X}, one has the {\it $\XX$-cluster variety} $\calX$ associated to $\boldsymbol{\Sigma}^\XX$, obtained by gluing the  
$\QQ$-rational tori $\{{\rm Spec}(\QQ[{\bf y}_t^{\pm 1}])\}_{t \in \TT_n}$ via the
birational morphisms $\{\phi_{t', t}^\XX\}_{t, t' \in \TT_n}$ 
as compositions of the one-step mutations in \eqref{eq:Y-mut}.
The algebra of regular functions on $\calX$ is  the upper cluster algebra
\[
\calU(\calX) = \bigcap_{t \in \TT_n} \QQ[{\bf y}_t^{\pm 1}].
\]
Elements in $\mathfrak{X}(\calX):=\{y_{t; i}: t \in \TT_n, \, i \in [1, n]\}$ are called
{\it $y$-variables} of $\calX$. 
Each ${\bf y}_t$, for $t \in \TT_n$, is called a cluster of $\calX$. The set of all clusters of $\calX$ is denoted by $\calC(\calX)$.
The assignment $\{B_t\}_{t \in \TT_n}$ is called the {\it matrix pattern} of $\calX$.
Unlike the case of an $\AA$-cluster variety $\calA$, 
where the celebrated {\it Laurent phenomenon} \cite{fz_2002}
says that each cluster variable of $\calA$ is a regular function on $\calA$, a $y$-variable of $\calX$ is not necessarily a regular function, i.e., we do not always 
have $\fX(\calX) \subset \calU(\calX)$. We thus introduce
\begin{equation}\label{eq:global-var}
\fX^{\rm global}(\calX) = \fX(\calX) \cap \calU(\calX),
\end{equation}
and we call elements $\fX^{\rm global}(\calX)$ {\it global $y$-variables} of $\calX$.

Suppose now that $\calA$ is an $\AA$-cluster variety associated to an
$\AA$-cluster pattern $\bSigma^\AA$ of rank $n$ in $\FF_{n+m}$, and $\calX$ is an $\XX$-cluster variety 
associated to 
an $\XX$-cluster pattern $\bSigma^\XX$ of rank $n$ in $\FF_n$, 
such that the matrix pattern $\{B_t\}_{t \in \TT_n}$ of $\calX$
is the principal part of the matrix pattern $\{\widetilde{B}_t\}_{t \in \TT_n}$ of $\calA$. 
Writing 
$\bSigma^\AA$  as in \eqref{eq:bold-Sigma-A}, the  {\it $\widehat{y}$-variables} of $\bSigma^\AA$ at $t\in \TT_n$ are defined as
\[
\widehat{\bf y}_t = {\bf x}_t^{B_t} {\bf p}^{P_t} \in \FF_{n+m}.
\]
Writing $\bSigma^\XX$  as in \eqref{eq:bold-Sigma-X}, 
it is well-known that one has the well-defined morphism $p: \calA \to \calX$ such that
\begin{equation}\label{eq:p-yt}
p^* ({\bf y}_t) = \widehat{\bf y}_t,  \hs t \in \TT_n.
\end{equation}
The triple $(\calA, \calX, p)$ is called a {\it cluster ensemble}, and $p$ is called the {\it ensemble map} 
\cite{FG:ensembles, GHK}. 

A cluster ensemble $(\calA, \calX, p)$ is said to have trivial coefficients
if $\calA$ does, i.e., if $m = 0$. For a cluster ensemble $(\calA, \calX, p)$ with trivial coefficients, if
$\calA$ is associated to the $\AA$-cluster pattern $\bSigma^\AA = \{({\bf x}_t, B_t)\}_{t \in \TT_n}$
and $\calX$ is associated to the $\XX$-cluster pattern 
$\bSigma^\XX = \{({\bf y}_t, B_t)\}_{t \in \TT_n}$, we write $\Sigma_t = 
({\bf x}_t, {\bf y}_t, B_t)$ for $t \in \TT_n$, and we call
$\bSigma = \{\Sigma_t = ({\bf x}_t, {\bf y}_t, B_t)\}_{t \in \TT_n}$
the {\it cluster pattern} of the ensemble $(\calA, \calX, p)$. 

While we will not be concerned with cluster ensembles with coefficients in general, we will need to 
consider {\it principal extensions} of cluster ensembles with trivial coefficients. More precisely, 
let $(\calA, \calX, p)$ be a cluster ensemble with trivial coefficients, and let 
\[
\bSigma = \{\Sigma_t = ({\bf x}_t = (x_{t; 1}, \ldots, x_{t; n}), \; {\bf y}_t
=(y_{t; 1}, \ldots, y_{t; n}), \;  B_t)\}_{t \in \TT_n}
\]
be its cluster pattern. For a given $t_0 \in \TT_n$, 
let $\{\widetilde{B}_t=(B_t\pp C_t)\}_{t \in \TT_n}$ be the matrix pattern such that
$\widetilde{B}_{t_0} = (B_{t_0} \pp I_n)$. Let $\Aprin$ be the $\AA$-cluster variety
associated to an 
$\AA$-cluster pattern
\begin{equation}\label{eq:Sigma-A-prin}
\bSigma^{\AA, {\rm prin}} = \{\Sigma_t^{\AA, {\rm prin}} =
(({\bf u}_t =(u_{t; 1}, \ldots, u_{t; n}), \, {\bf p} = (p_1, \ldots, p_n)), \; \widetilde{B}_t)\}_{t \in \TT_n}
\end{equation}
of rank $n$ in a field extension $\FF_{2n}$ of $\QQ$ of pure transcendental degree $2n$. One then has 
the  embedding $\iota: \calA \to \Aprin$ and the projection
$p^{\rm prin}: \calA^{\rm prin} \to \calX$,
such that 
\begin{equation}\label{eq:iota-p-prin}
\iota^* ({\bf u}_{t}) = {\bf x}_{t} \hs \mbox{and} \hs 
(p^{\rm prin})^*({\bf y}_t) ={\bf u}_t^{B_t}{\bf p}^{C_t},
\hs t \in \TT_n.
\end{equation}
Note that $p = p^{\rm prin}\circ \iota: \calA \to \calX$,
and $\iota(\calA)=
\{a \in\Aprin: p_1(a) = \cdots =p_n(a) = 1\}$. Let 
$(\QQ^\times)^n$ act on $\calU(\Aprin)$ such that 
${\bf u}_{t_0; i}$ has weight $e_i$ and $p_i$ has weight $-B_{t_0}e_i$ for $i \in [1, n]$.
Then \cite[Lemma 7.10]{GHKK}
\[
(p^{\rm prin})^*: \;\; \calU(\calX) \longrightarrow  \calU(\Aprin)
\]
identifies $\calU(\calX)$ with the sub-algebra of $(\QQ^\times)^n$-invariant elements of $\calU(\Aprin)$.

\bde{:prin-extension}
{\rm With the notation as above, we call $\Aprin$ a {\it principal extension of $\calA$},
and $(\Aprin, \calX, p^{\rm prin})$ a {\it principle extension} of $(\calA, \calX, p)$ at
the seed of $\calA$ at $t_0$.
\hfill $diamond$
}
\ede

\subsection{DT transformations and characteristic polynomials at positive fixed points}\label{ss:DT-expon}
Let $(\calA, \calX, p)$ be a rank $n$ cluster ensemble with trivial coefficients, and let
$\{B_t\}_{t \in \TT_n}$ be its matrix pattern. 
Let
$I_n$ be the $n \times n$ identity matrix. 
Recall from $\S$\ref{ss:nota} that for $\sigma \in S_n$ we have 
$\overline{\sigma}\in {\rm GL}(n, \ZZ)$.

\bde{:reddening} \cite{keller_2012, keller_demonet_2020}
The matrix pattern $\{B_t\}_{t \in \TT_n}$ of $(\calA, \calX, p)$ is said to 
{\it admit reddening sequences}
if there exist $t_0, t_0[1] \in \TT_n$ such that the unique matrix pattern 
$\{\widetilde{B}_t = (B_t\pp C_t)\}_{t \in \TT_n}$
with $\widetilde{B}_{t_0} = (B_{t_0}\pp I_n)$ has the property that 
$\widetilde{B}_{t_0[1]}= 
(\overline{\sigma}\, B_{t_0} \overline{\sigma}^{\, -1}\pp -\overline{\sigma}^{\, -1})$ for 
some $\sigma$ in the permutation group $S_n$. 
\ede

Given $\sigma \in S_n$, let 
$\TT_n^\sigma$ be the $n$-regular tree with the same set of vertices as $\TT_n$ but the edges of $\TT_n^\sigma$ are labeled in such a way that for
$t, v \in \TT_n$ and $k \in [1, n]$,
\[
t \,\stackrel{k}{\mathdash\!\mathdash}\, v \;\; \mbox{in} \;\; \TT_n^\sigma \hs \mbox{if and only if}\hs
t \,\stackrel{\sigma(k)}{\mathdash\!\mathdash}\, v\;\; \mbox{in} \;\;\TT_n.
\]
Given any $t_0, t_0[1] \in \TT_n$ and $\sigma \in S_n$, one then has a unique tree isomorphism 
from $\TT_n$ to $\TT_n^\sigma$ that 
sends $t_0$ to $t_0[1]$ and we denote it by $\TT_n \to \TT_n^\sigma: t \mapsto t[1]$. More precisely, if
$t_0 \,\stackrel{i_1}{\mathdash\!\mathdash}  \bullet\;\;\cdots \;\;\bullet\,
\stackrel{i_{p}}{\mathdash\!\mathdash}{t}$ in $\TT_n$, then 
\[
t_0[1] \,\stackrel{i_1}{\mathdash\!\mathdash}  \bullet\;\;\cdots \;\;\bullet\,
\stackrel{i_{p}}{\mathdash\!\mathdash}{t[1]} \;\; \mbox{in}\; \;\TT_n^\sigma, \hs \mbox{equivalently},\hs 
t_0[1] \,\stackrel{\sigma(i_1)}{\mathdash\!\mathdash}  \bullet\;\;\cdots \;\;\bullet\,
\stackrel{\sigma(i_{p})}{\mathdash\!\mathdash}{t[1]} \;\; \mbox{in} \; \; \TT_n.
\]

\bde{:reddening-triple}
{\rm 
In the setting of \deref{:reddening}, we
call $(t_0, t_0[1], \sigma)$ a {\it reddening triple} of the matrix pattern $\{B_t\}_{t \in \TT_n}$, 
or of the cluster ensemble $(\calA, \calX, p)$.
\hfill $\diamond$
}
\ede

\bre{:triple-not-unique}
{\rm
 The permutation $\sigma$ in a reddening triple $(t_0, t_0[1], \sigma)$ is uniquely determined by 
 $(t_0, t_0[1])$, as $\overline{\sigma}^{\, -1} = -C_{t_0[1]}$.  On the other hand,
reddening triples for $(\calA, \calX, p)$, when they exist, are not unique. Indeed, given $t_0 \in \TT_n$
there could be more than one pair $(t_0[1], \sigma)$ such that $(t_0, t_0[1], \sigma)$ is a reddening triple. 
Moreover, if $(t_0, t_0[1], \sigma)$ is a reddening
triple, so is $(t, t[1], \sigma)$ for every $t \in \TT_n$, where $\TT_n \to \TT_n^\sigma, t\to t[1],$ 
is the tree isomorphism sending $t_0$ to $t_0[1]$. In particular, for every $t \in \TT_n$ one has 
\begin{equation}\label{eq:Btt}
B_{t[1]} = \overline{\sigma} B_t \overline{\sigma}^{\, -1}, \hs t \in \TT_n.
\end{equation}
Existence of reddening sequences of the matrix pattern $\{B_t\}_{t \in \TT_n}$ can also be characterized 
by the existence of a (unique) sink in its oriented exchange graph, and we refer to 
\cite[\S 2]{Brustle-Maximal-green} and \cite{Muller-existence-green-seq} for details.
\hfill $\diamond$
}
\ere

For the following statement and definition on DT transformations, we refer to 
\cite[Proposition 2.3.3 and Corollary 2.3.4]{Qin:22},  \cite[Theorem 3.2]{Goncharov-Shen-DT}, 
\cite[Theorem 6.5, Remark 6.6]{keller_2012}, and \cite[\S 3.3]{keller_demonet_2020}.

\bth{:DT} Let $(\calA, \calX, p)$ be a cluster ensemble with trivial coefficients and cluster pattern
\begin{equation}\label{eq:bSigma}
\bSigma = \{\Sigma_t = ({\bf x}_t = (x_{t; 1}, \ldots, x_{t; n}), \; {\bf y}_t
=(y_{t; 1}, \ldots, y_{t; n}), \;  B_t)\}_{t \in \TT_n},
\end{equation}
and assume that $\{B_t\}_{t \in \TT_n}$ admits reddening sequences. Given any reddening triple
$(t_0, t_0[1], \sigma)$ of $\{B_t\}_{t \in \TT_n}$,
there are biregular automorphisms
$\DTA: \calA \to \calA$ and $\DTX: \calX \to \calX$, uniquely determined by
\begin{align}\label{eq:DTA-x}
(\DTA)^* ({\bf x}_t) &= {\bf x}_{t[1]} \,\overline{\sigma} = (x_{t[1]; \sigma(1)}, \, \ldots, \, 
x_{t[1]; \sigma(n)}),\\
\label{eq:DTX-y}
(\DTX)^* ({\bf y}_t) &= {\bf y}_{t[1]}\, \overline{\sigma} =
 (y_{t[1]; \sigma(1)}, \, \ldots, \, y_{t[1]; \sigma(n)})
\end{align}
for every $t \in \TT_n$, where  $t \mapsto t[1]$ is the 
unique tree isomorphism  $\TT_n \to \TT_n^\sigma$ sending $t_0$ to $t_0[1]$.
The two maps $\DTA$ and $\DTX$ are independent of the choice of the reddening triple $(t_0, t_0[1], \sigma)$.
\eth

\bde{de:DT}
When the matrix pattern of the 
cluster ensemble $(\calA, \calX, p)$ (with trivial coefficients) admits reddening sequences, 
we also say that $(\calA, \calX, p)$ admits DT transformations, and the 
two maps 
\[
\DTA: \;\;\calA \longrightarrow \calA \hs \mbox{and} \hs 
\DTX: \;\; \calX\longrightarrow \calX
\]
are respectively called the {\it DT transformation of $\calA$} and 
the {\it DT transformation of
$\calX$}. 
\ede

For the rest of $\S$\ref{ss:DT-expon}, we fix a cluster ensemble $(\calA, \calX, p)$ 
that admits DT transformations, and let the cluster pattern of $(\calA, \calX, p)$ be as in \eqref{eq:bSigma}.


\ble{:DT-p}
One has $p \circ \DTA = \DTX \circ p: \calA \to \calX$.
\ele

\begin{proof}
Let $t \in \TT_n$ be arbitrary. In the notation of \eqref{eq:DTA-x} and \eqref{eq:DTX-y} and by
\eqref{eq:a-VW}, \eqref{eq:a-sigma} and \eqref{eq:Btt}, 
\begin{align*}
\DTA^* (p^*({\bf y}_{t})) &= \DTA^*(\widehat{\bf y}_{t}) = \DTA^*({\bf x}_{t}^{B_{t}})
=({\bf x}_{t[1]}\overline{\sigma})^{B_{t}} 
= {\bf x}_{t[1]}^{\overline{\sigma}{B_{t}}} = 
{\bf x}_{t[1]}^{B_{t[1]}\overline{\sigma}} =
\left({\bf x}_{t[1]}^{B_{t[1]}}\right){\overline{\sigma}} \\
&= \widehat{\bf y}_{t[1]} \overline{\sigma} = p^*({\bf y}_{t[1]}) \overline{\sigma}
= p^*({\bf y}_{t[1]} \overline{\sigma}) = p^* (\DTX^*({\bf y}_{t})).
\end{align*}
\end{proof}

Consider now a principal extension 
$(\Aprin, \calX, p^{\rm prin})$  of $(\calA, \calX, p)$ at the seed of $\calA$ at some $t_0 \in \TT_n$
(see \deref{:prin-extension}). Write again the 
$\AA$-cluster pattern $\bSigma^{\AA, {\rm prin}}$ in $\FF_{2n}$ as  
\begin{equation}\label{eq:Sigma-A-prin-1}
\bSigma^{\AA, {\rm prin}} = \{\Sigma_t^{\AA, {\rm prin}} =
(({\bf u}_t =(u_{t; 1}, \ldots, u_{t; n}), \, {\bf p} = (p_1, \ldots, p_n)), \; \widetilde{B}_t)\}_{t \in \TT_n},
\end{equation}
where $\widetilde{B}_{t_0} = (B_{t_0} \pp I_n)$. Recall from \eqref{eq:iota-p-prin}
the embedding $\iota: \calA \to \Aprin$
and the projection $p^{\rm prin}: \Aprin \to \calX$.

\bld{:psi-extension}
{\rm
Given $(t_0[1], \sigma)$ such that $(t_0, t_0[1], \sigma)$ is a reddening triple of $(\calA, \calX, p)$,
for any $n \times n$ integer matrices $R$ and $E$ satisfying 
\begin{equation}\label{eq:RE}
RB_{t_0} + E = -I_n,
\end{equation}
there is a bi-regular automorphism $\psi: \Aprin \to \Aprin$ uniquely determined by 
\begin{equation}\label{eq:psi-A-prin}
\psi^*({\bf u}_{t_0}) = 
({\bf u}_{t_0[1]}\overline{\sigma}){\bf p}^R \hs \mbox{and}\hs
\psi^*({\bf p}) = {\bf p}^E.
\end{equation}
Moreover, all three squares in the following diagram commute:
\begin{equation}\label{eq:diagrams}
\begin{tikzcd}
\calA  \arrow[r, "\iota" ] \arrow[d,"\DTA" '] &\calA^{\rm prin}\arrow[d, "\psi"] \arrow[r, "p^{\rm prin}"] & \calX \arrow[d,"\DTX"] \\
\calA  \arrow[r, "\iota" '] & \Aprin \arrow[r,"p^{\rm prin}" '] & \calX
\end{tikzcd}.
\end{equation}
Following \cite[\S 5.1]{kimura-qin-wei}, we call $\psi: \Aprin\to \Aprin$
a {\it twisted automorphism of $\Aprin$ of DT-type}.
\hfill $\diamond$
}
\eld

\begin{proof}
Write $\widetilde{B}_t = (B_t\pp C_t)$  for $t \in \TT_r$, and consider the  $\AA$-cluster pattern
\[
(\bSigma^{\AA, {\rm prin}})^\sigma = 
\left\{(\Sigma_t^{\AA, {\rm prin}})^\sigma :=(({\bf u}_t^\sigma, \; 
{\bf p}), \; \widetilde{B}^\sigma_t)\right\}_{t \in \TT_n^\sigma}
\]
in $\FF_{2n}$,
where ${\bf u}_t^\sigma ={\bf u}_t \overline{\sigma}$ and
$\widetilde{B}^\sigma_t=(\overline{\sigma}^{\,-1} B_t \overline{\sigma} \pp (C_t \overline{\sigma}))$ 
for $t \in \TT_n^\sigma$. It is evident that $(\bSigma^{\AA, {\rm prin}})^\sigma$ gives rise to the same 
$\AA$-cluster variety $\calA^{\rm prin}$. 
By \eqref{eq:a-VW} -\eqref{eq:ab-V-sigma}, the $\widehat{y}$-variables of 
$(\bSigma^{\AA, {\rm prin}})^\sigma$ at $t \in \TT_n^\sigma$ are
\[
\widehat{\bf y}^\sigma_t = ({\bf u}_t \overline{\sigma})^{\overline{\sigma}^{\, -1} B_t \overline{\sigma}}
{\bf p}^{C_t \overline{\sigma}} 
= {\bf u}_t^{B_t \overline{\sigma}}
{\bf p}^{C_t \overline{\sigma}} 
= ({\bf u}_t^{B_t}{\bf p}^{C_t}){\overline{\sigma}} = 
\widehat{\bf y}_t \overline{\sigma}.
\]
Let $\Psi: \FF_{2n} \to \FF_{2n}$ be the unique field isomorphism defined by
\begin{equation}\label{eq:Psi-A-prin}
\Psi({\bf u}_{t_0}) = {\bf u}_{t_0[1]}^\sigma{\bf p}^R = 
({\bf u}_{t_0[1]}\overline{\sigma}){\bf p}^R, \hs 
\psi^*({\bf p}) = {\bf p}^E.
\end{equation}
As $C_{t_0} = -\overline{\sigma}^{\, -1}$, one has
$\widehat{\bf y}^\sigma_{t_0[1]} = {\bf u}_{t_0[1]}^{B_{t_0[1]}\overline{\sigma}} {\bf p}^{-1}$. 
One the other hand,
\[
\Psi(\widehat{\bf y}_{t_0}) = \Psi({\bf u}_{t_0}^{B_{t_0}} {\bf p}) = 
(({\bf u}_{t_0[1]}\overline{\sigma}) {\bf p}^R)^{B_{t_0}} {\bf p}^{E} =
({\bf u}_{t_0[1]}\overline{\sigma})^{B_{t_0}} {\bf p}^{RB_{t_0}+E} =    
{\bf u}_{t_0[1]}^{B_{t_0[1]}\overline{\sigma}} {\bf p}^{RB_{t_0}+E}.
\]
Condition \eqref{eq:RE} is thus equivalent to 
$\Psi(\widehat{\bf y}_{t_0}) = \widehat{\bf y}^\sigma_{t_0[1]}$.
By \cite{Fraser:quasi} (see also \cite[$\S$4.3]{GHL:friezes} and \cite[\S 5.1]{kimura-qin-wei}),
 $\Psi: \FF_{2n} \to \FF_{2n}$ restricts to
an algebra automorphism of $\calU(\calA^{\rm prin})$, still denoted as $\Psi$,
which is 
an {\it upper cluster algebra quasi-automorphism} of $\calU(\calA^{\rm prin})$ in the sense of 
\cite{Fraser:quasi}. In particular, for each $t \in T$ there exists an $n \times n$ integer matrix $R_t$
such that
\begin{equation}\label{eq:Psi}
\Psi({\bf u}_t) = {\bf u}_{t[1]}^{\overline{\sigma}}{\bf p}^{R_t}, \hs \Psi(\widehat{\bf y}_{t}) = \widehat{\bf y}^\sigma_{t[1]}=\wh{\bf y}_{t[1]}\overline{\sigma}, \hs \hs t \in \TT_r,
\end{equation}
where again $\TT_n \to \TT_n^\sigma, t \to t[1]$, is the tree isomorphism sending $t_0$ to $t_0[1]$.
Let $\psi:\calA^{\rm prin} \rightarrow \calA^{\rm prin}$
be the automorphism of the variety $\calA^{\rm prin}$ such that $\psi^* = \Psi: \calU(\Aprin) \to \calU(\Aprin)$. Then  \eqref{eq:psi-A-prin} holds, and  it follows from \eqref{eq:Psi}, \leref{:DT-p} and the identity $p = p^{\rm prin} \circ \iota$ that the diagrams in \eqref{eq:diagrams} commute.
\end{proof}

Recall now that the set of {\it positive points} of $\calA$ and of $\calX$ are defined to be
\[
    \calApos = \{ {a} \in \calA(\mathbb{R}): \;x({ a}) >0, \forall x \in \fX(\calA)\},\hs 
    \calXpos = \{ {b} \in \calX(\mathbb{R}): \;y({b}) >0, \forall y \in \fX(\calX)\}.
\]
Each cluster ${\bf x} \in \calC(\calA)$, respectively ${\bf y} \in \calC(\calX)$,
gives rise to global coordinates 
\begin{align}\label{eq:pos-coords}
    {\bf x}:\; \calApos \longrightarrow (\RR_{>0})^n, \quad {a} \mapsto {\bf x}({a}), \quad \text{respectively} \quad {\bf y}:\; \calXpos \longrightarrow (\RR_{>0})^n, \quad { b} \mapsto {\bf y}({ b}),
\end{align}
and endows $\calApos$, respectively $\calXpos$,  with the unique smooth structure for which the cluster variables of $\calA$, respectively the $y$-variables of $\calX$, are smooth functions. 
By \eqref{eq:DTA-x} and \eqref{eq:DTX-y},  one has the diffeomorphisms
\[
    \DTApos:=\DTA|_{\calA_{>0}}: \;\calApos \longrightarrow \calApos \quad \text{ and } \quad 
    \DTXpos:=\DTX|_{\calX_{>0}}:\;  \calXpos \longrightarrow \calXpos.
\]
By \leref{:DT-p}, if ${\bf a} \in \calApos$ is a fixed point of $\DTApos$, then 
${\bf b} = p({\bf a}) \in \calXpos$ is a fixed point of $\DTXpos$. 
Taking the linearizations of $\DTApos$ and $\DTXpos$ at ${\bf a}$ and ${\bf b}$ respectively, we have
the real linear maps
\begin{align*}
    &\dtAa :={\rm d}_{\bf a}(\DTApos): \;{\rm T}_{\bf a} (\calApos) \longrightarrow {\rm T}_{\bf a} (\calApos), \\
    &\dtXb :={\rm d}_{\bf b}(\DTXpos):\; {\rm T}_{\bf b} (\calXpos) \longrightarrow {\rm T}_{\bf b} (\calXpos).
\end{align*}
Let now $\Aprin$ be a principal extension of $\calA$, and consider also the positive part 
\[
\Aprinpos = \{a \in \Aprin: \, u(a) > 0 \; \forall \; u \in \fX(\Aprin), \, p_1(a) > 0, \ldots, p_n(a) > 0\}
\]
of $\calA^{\rm prin}$.  
Let $\psi: \Aprin \to \Aprin$ be a twisted automorphism of $\Aprin$ of DT-type. By 
\eqref{eq:psi-A-prin}, one has 
$\psi (\Aprinpos) \subset \Aprinpos$. Denote the restriction of $\psi$ to $\Aprinpos$  by
$\psi_{>0}: \Aprinpos \rightarrow \Aprinpos$.
Identifying  $\calA \cong \iota(\calA)\subset \Aprin$, then 
$\calApos \subset \Aprinpos$. In particular, ${\bf a} \in \Aprinpos$ is fixed by $\psipos$, and 
 we have the linear map
\[
d_{\bf a}(\psi_{>0}): \;\; T_{\bf a} (\Aprinpos) \longrightarrow T_{\bf a} (\Aprinpos). 
\]

\bth{:same-char} Let $(\calA, \calX, p)$ be a rank $n$ cluster ensemble (with trivial coefficients) admitting 
DT transformations, let ${\bf a} \in \calApos$ be any fixed point of $\DTApos$, and let
${\bf b} = p({\bf a}) \in \calXpos$. Let $\Aprin$ be any principal extension of $\calA$, 
and let $\psi: \Aprin \to \Aprin$ be any twisted automorphism of $\Aprin$ of DT-type. Let $E$ be the 
$n \times n$ integer matrix such that 
$\psi^*({\bf p}) = {\bf p}^E$, where ${\bf p} = (p_1, \ldots, p_n)$ are the frozen variables of $\Aprin$. 

1) The two linear maps $\dtAa$ and $\dtXb$ have the same characteristic polynomial 
$\calP(x) \in \RR[x]$;

2) The characteristic polynomial of $d_{\bf a}(\psi_{>0})$ is equal to $\calP(x) \calP_E(x)$,
where $\calP_E(x) = \det(x-E)$.
\eth

\begin{proof} We compute the dual map 
$d_{\bf a}^*(\psi_{>0}):  T_{\bf a}^* (\Aprinpos) \to T_{\bf a}^* (\Aprinpos)$ of 
$d_{\bf a}(\psi_{>0})$ in
two bases of the cotangent space $T_{\bf a}^* (\Aprinpos)$. 
To this end, let $(t_0, t_0[1], \sigma)$ be a reddening triple of $(\calA, \calX, p)$  such that
$\Aprin$ is a principal extension of $\calA$ at
$t_0$ and that \eqref{eq:psi-A-prin} holds for $R$ and $E$ satisfying \eqref{eq:RE}. 
Let the $\AA$-cluster pattern of $\Aprin$ be as in \eqref{eq:Sigma-A-prin-1}. Setting
\[
{\bf u} = {\bf u}_{t_0} = (u_1, \ldots, u_n) \hs \mbox{and} \hs
\widehat{\bf y} = 
\widehat{\bf y}_{t_0} = (\widehat{y}_1, \ldots, \widehat{y}_n),
\]
we will regard $({\bf u}, {\bf p})$ and $({\bf u}, \wh{\bf y})$ as two coordinate systems on $\Aprinpos$
via the parametrizations
\begin{align*}
&\Aprinpos \ni a \longmapsto ({\bf u}(a), {\bf p}(a)) = (u_1(a), \ldots, u_n(a),\, p_1(a), \ldots, p_n(a)) \in (\RR_{>0})^{2n}, \\
&\Aprinpos \ni a \longmapsto ({\bf u}(a), \wh{\bf y}(a)) = (u_1(a), \ldots, u_n(a), \,\widehat{y}_1(a), \ldots, \widehat{y}_n(a)) \in (\RR_{>0})^{2n}.
\end{align*}
We thus have two bases ${\bf w}_1$ and ${\bf w}_2$ of $T_{\bf a}^* (\Aprinpos)$, respectively given by
\begin{align*}
{\bf w}_1  &= (d_{\bf a} {\bf u} = (d_{\bf a}u_1, \ldots, d_{\bf a}u_n), \,  d_{\bf a} p_1, \ldots, 
d_{\bf a}p_n),\\
{\bf w}_2 & = ({\bf u}({\bf a})^{-1}d_{\bf a} {\bf u}, \; d_{\bf a} \wh{\bf y}) 
=\left(\frac{d_{\bf a}u_1}{u_1({\bf a})}, \ldots, \frac{d_{\bf a}u_n}{u_n({\bf a})}, \;  d_{\bf a} \wh{y}_1, \ldots, d_{\bf a}\wh{y}_n\right).
\end{align*}
We also regard $d_{\bf a}{\bf u}$ as a basis of $T_{\bf a}^* (\calApos)$,
and let $M$ be the $n \times n$ real matrix such that
\[
\dtAa^* (d_{\bf a}u_1, \ldots, d_{\bf a}u_n) = (d_{\bf a}u_1, \ldots, d_{\bf a}u_n) M.
\]
Write ${\bf u}_{t_0[1]}\overline{\sigma} = (U_1({\bf u}, {\bf p}), \ldots, U_n({\bf u}, {\bf p}))$,
where $U_i({\bf u}, {\bf p}) \in \ZZ[{\bf u}^{\pm 1}, {\bf p}]$ for each $i \in [1, n]$. 
Noting that $\psipos|_{\calApos} = \DTApos$ and $p_1({\bf a}) = \cdots = p_n({\bf a}) = 1$, 
it then follows from \eqref{eq:psi-A-prin} and \eqref{eq:d-xV-2} that 
\begin{equation}\label{eq:d-psi-w1}
d_{\bf a}^* (\psipos) {\bf w}_1 = {\bf w}_1 \left(\begin{array}{cc} M & 0 \\ M' & E\end{array}\right)
\end{equation}
for some $n \times n$ real matrix $M'$. To compute $d_{\bf a}^* (\psipos) {\bf w}_2$, we use the 
Separation Formula \cite[Corollary 6.3]{FZ:ClusterIV} to write 
${\bf u}_{t_0[1]} = {\bf u}_{t_0}^G F(\wh{\bf y})$, where 
$G \in {\rm  GL}(n, \ZZ)$ and $F = (F_1, \ldots, F_n)$ are respectively the $g$-vector matrix 
and the $F$-polynomials of $\{\widetilde{B}_t\}_{t \in \TT_n}$ at $t_0[1]$ with respect to $t_0$. 
As $C_{t_0[1]} = -\overline{\sigma}^{\, -1}$,
by tropical duality \cite[(3.11)]{NZ12} (see also \cite[Proposition 2.3.3]{Qin:22}), 
we have $G = -\overline{\sigma}^{\, -1}$.  Thus 
\[
(\psipos)^*({\bf u})
= \left(({\bf u}_{t_0}^{-\overline{\sigma}^{\, -1}} F(\wh{\bf y}))\overline{\sigma}\right){\bf p}^R = 
({\bf u}_{t_0}^{-\overline{\sigma}^{\, -1}} \overline{\sigma})(F(\wh{\bf y})\overline{\sigma})({\bf u}^{-B_{t_0}}\wh{\bf y})^R 
= {\bf u}_{t_0}^{-(I_n+B_{t_0}R)} \wh{\bf y}^R(F(\wh{\bf y})\overline{\sigma}).
\]
In the $({\bf u}, \wh{\bf y})$-coordinates, we then have
\begin{equation}\label{eq:psi-u-hty}
(\psipos)^*({\bf u}, \wh{\bf y})
 = ({\bf u}_{t_0}^{-(I_n+B_{t_0}R)} \wh{\bf y}^R(F(\wh{\bf y})\overline{\sigma}), \; (\psipos)^*(\wh{\bf y})).
\end{equation}
Consider now the cluster  ${\bf y}_{t_0}$ of $\calX$ and write $(y_1, \ldots, y_n) = {\bf y}_{t_0}|_{\calXpos}$ 
as a coordinate system on $\calXpos$, so that $(d_{\bf b}y_1, \ldots,
d_{\bf b}y_n)$ is a basis of $T_{\bf b} (\calXpos)$. Let $N$ be the real $n\times n$ matrix such that 
\[
\dtXb^* (d_{\bf b}y_1, \ldots,d_{\bf b}y_n) = (d_{\bf b}y_1, \ldots,d_{\bf b}y_n) N.
\]
It follows from \eqref{eq:d-xV-2}, \eqref{eq:diagrams} and \eqref{eq:psi-u-hty} that 
\begin{equation}\label{eq:d-psi-w2}
d_{\bf a}^* (\psipos) {\bf w}_2 = {\bf w}_2 \left(\begin{array}{cc} -(I_n+B_{t_0}R) & 0 \\ N' & N\end{array}\right)
\end{equation}
for some real $n \times n$ matrix $N'$. 
Denote the characteristic polynomials of $\dtAa$, $\dtXb$ and $d_{\bf a}(\psi_{>0})$ respectively 
by $\calP_{\calA}(x)$, $\calP_{\calX}(x)$ and 
$\calP_{\psipos}(x)$, and let $\calP_E(x)$ and $\calQ(x)$ be the respective  
characteristic polynomials of $E$ and $-I_n-B_{t_0}R$. 
 It now follows from
\eqref{eq:d-psi-w1} and \eqref{eq:d-psi-w2} that 
\begin{equation}\label{eq:PPQ}
\calP_{\psipos}(x) = \calP_E(x) \calP_{\calA}(x) = \calQ(x) \calP_{\calX}(x).
\end{equation}
By \eqref{eq:RE} we have $E = -I_n-RB_{t_0}$. Thus
\[
\calP_E(x) = \det(xI_n-E) = \det((x+1)I_n+ RB_{t_0}) \hs \mbox{and} \hs
\calQ(x) = \det((x+1)I_n + B_{t_0}R).
\]
As $\det(xI_n+AB) = \det(xI_n+BA)$  for any two $n \times n$ matrices $A$ and $B$, we 
have $\calP_E(x) = \calQ(x)$. We thus conclude from \eqref{eq:PPQ} that $\calP_{\calA}(x) = \calP_{\calX}(x)$ and 
$\calP_{\psipos}(x) =  \calP_{\calA}(x)\calP_E(x)$.
\end{proof}

As in \deref{:cluster-expo-intro}, when 
$\DTA$  has a unique fixed point ${\bf a} \in \calApos$ and has period ${\bf q}$,
we call the exponents of $\dtAa$ with respect to ${\bf q}$ the {\it cluster exponents} of
$(\calA, \calX, p)$ with respect to ${\bf q}$.

\ble{:fixpt-frozen-1}
Let $\psi: \Aprin \to \Aprin$ be a twisted automorphism of $\Aprin$ of DT-type, with $\psi^*({\bf p}) = {\bf p}^E$ for some $E \in M_n(\ZZ)$, and suppose that $\psi({\bf a}) = {\bf a} \in \Aprinpos$. If $E$ does not have eigenvalue 1, then ${\bf p}({\bf a}) =(1,\ldots , 1)$. 
\ele
\begin{proof}
    If $\psi({\bf a}) = {\bf a}$, we have ${\bf p}({\bf a}) = {\bf p}({\bf a})^E$, and thus $\log ({\bf p}({\bf a})) = \log ({\bf p}({\bf a})) E$, from which the claim is immediate.
\end{proof}

\section{Fixed points of DT transformations and cluster exponents in finite type}\label{s:fixedpoint}

\subsection{DT transformations on cluster ensembles of finite type}\label{ss:finitetype} 
Recall from \cite{FZ:II} that the  {\it Cartan companion} of an $n\times n$
mutation matrix $B = (b_{i,j})$ is the matrix $(a_{i,j})$ given by 
\[
 a_{i,i} = 2  \;\; \forall i \in [1, n] \hs \mbox{and} \hs    a_{i,j} = - |b_{i,j}| \;\;\forall i \neq j.
\]
In \cite{FZ:II}, Fomin and Zelevinsky showed that a cluster ensemble $(\calA, \calX, p)$ (with trivial coefficients)
with matrix
pattern $\{B_t\}_{t \in \TT_n}$ is 
of finite type, i.e., $\calA$, equivalently $\calX$, has finitely many clusters, if and only 
if there exists a $t \in \TT_n$ such that the Cartan companion $A_t$ 
of $B_t$ is a Cartan matrix of finite type, and that the Cartan-Killing type of $A_t$ is
independent of such a $t \in \TT_n$.  In this way, each cluster ensemble 
$(\calA, \calX, p)$ of finite type 
(and with trivial coefficients) is associated to a root system $\Phi$ of unique Cartan-Killing type, and we
say that $(\calA, \calX, p)$ is of type $\Phi$.

Throughout the rest of $\S$\ref{s:fixedpoint}, 
we fix a finite (reduced, crystallographic) indecomposable root system $\Phi$ of rank $n$ and 
a choice of simple roots in $\Phi$. Denote the corresponding Cartan matrix by
$A = (a_{i, j})_{i, j \in [1, n]}$, 
and introduce the mutation matrix $B_A = (b_{i,j})_{i, j\in [1, n]}$, where for $i, j \in [1, n]$,
\begin{equation}\label{eq:BA}
    b_{i,j} = \begin{cases}
        0 & \text{ if } i = j, \\
        -a_{i,j} & \text{ if } i>j, \\
        a_{i,j} & \text{ if } i< j.
    \end{cases}
\end{equation}
In the notation in \eqref{eq:LU-intro} and with ${\bf 1} = (1, \ldots, 1)$, we have
$A = L_{\bf 1} + U_{\bf 1}$ and $B_A = -L_{\bf 1} + U_{\bf 1}$.

Let $\FF_n$ be a field extension of $\QQ$ of pure transcendental degree $n$, and let
\[
{\bf x} = (x_1, \ldots, x_n) \hs \mbox{and} \hs 
{\bf y} = (y_1, \ldots, y_n)
\]
be two ordered sets of free transcendental generators of $\FF_n$ over $\QQ$. 
Fixing $t_0 \in \TT_n$, we  then have the cluster ensemble $(\calA, \calX, p)$ with trivial coefficients
associated to the cluster pattern
\[
\bSigma_A=  
\{\Sigma_t = ({\bf x}_t=(x_{t; 1}, \ldots, x_{t; n}), \; {\bf y}_t=(y_{t; 1}, \ldots, y_{t; n}),\; B_t)\}_{t \in \TT_n},
\]
where $\Sigma_{t_0} = ({\bf x}, {\bf y}, B_A)$.  As $A$ is the Cartan companion of $B_A$, 
the cluster ensemble  $(\calA,\calX,p)$ is of type $\Phi$.
Up to isomorphisms, every cluster ensemble  of type $\Phi$ is obtained this way.

A direct computation gives
    \[
        \mu_n \cdots \mu_2 \mu_1 \begin{pmatrix}
            B_A \\ I_n
        \end{pmatrix} = \begin{pmatrix}
            B_A \\ -I_n
        \end{pmatrix}.
    \]
Thus $(\calA, \calX, p)$ admits DT transformations
$\DTA: \calA \to \calA$ and $\DTX: \calX \to \calX$. Let $t_0[1] \in \TT_n$ be given by the path 
$t_0 \,\stackrel{1}{\mathdash} t_1 \,\stackrel{2}{\mathdash} \cdots \, \stackrel{n-1}{\mathdash}
t_{n-1} \stackrel{n}{\mathdash} t_0[1]$ in $\TT_n$ and write 
    \[
    \Sigma_{t_0[1]} = \mu_n \cdots \mu_2 \mu_1(\Sigma_{t_0}) =
    ({\bf x}[1] = (x_{n+1}, \ldots, x_{2n}), \,\;{\bf y}[1] = (y_{n+1}, \ldots, y_{2n}),\;\, B_A).
    \]
By \thref{:DT}, the DT transformations $\DTA: \calA \to \calA$ and $\DTX: \calX \to \calX$
are uniquely determined by 
\[
 \DTA^*({\bf x}) = {\bf x}[1] \hs \mbox{and} \hs \DTX^*({\bf y}) = {\bf y}[1].
\]
The following periodicity property of $\DTA$ and $\DTX$ are a special case of the 
celebrated {\it Zamolodchikov periodicity} (conjectured in \cite{Zamolo:Y-system}, 
proved for finite types in \cite{FZ:associahedra} and in full generality in \cite{Keller:periodicity}). 

\ble{:DT-periodic}\cite{FZ:associahedra}
With ${\bf h}$ as the Coxeter number of $\Phi$, one has
    \[
    (\DTA)^{{\bf h}+2} = (\DTX)^{{\bf h}+2} = {\rm Id}.
    \]
When the longest element $w_0$ of the Weyl group of $\Phi$ acts as minus the identity map on $\Phi$, one has $(\DTA)^{\frac{{\bf h}+2}{2}} = (\DTX)^{\frac{{\bf h}+2}{2}} = {\rm Id}$.
\ele

To further understand $\DTA$ and $\DTX$, we first note that a direct computation using \eqref{eq:X-mut} gives
\begin{equation}\label{eq:x-Cox-mutation}
x_i x_{n+i} = 1 + \prod_{j = i+1}^n x_j^{-a_{j,i}} \prod_{j=1}^{i-1} x_{n+j}^{-a_{j,i}}.
\end{equation}
The cluster variables $x_1, \ldots, x_n, x_{n+1}, \ldots, x_{2n}$ of $\calA$ and the relations in \eqref{eq:x-Cox-mutation}
are part of the so-called {\it generic frieze pattern}  of $\calA$, which is the unique
surjective map \cite{CdsgL:ringel-conj, GHL:friezes}
\begin{equation}\label{eq:gen-frieze-pat}
[1, n] \times \ZZ \longrightarrow  \fX(\calA), \; \;
(i, m) \longmapsto x_{i, m},
\end{equation}
satisfying $x_{i, 0} = x_i$ for every $i \in [1, n]$, and 
\begin{equation}\label{eq:frieze-mut-0}
x_{i,m} \, x_{i,m+1} = 1 + \prod_{j=i+1}^n x_{j,m}^{-a_{j,i}} \, \prod_{j=1}^{i-1} x_{j,m+1}^{-a_{j,i}}
\end{equation}
for all $i \in [1, n]$ and $m \in \ZZ$. By \eqref{eq:x-Cox-mutation}, $x_{n+i} = x_{i, 1}$ for every $i \in [1, n]$.
An inductive argument then gives
\begin{equation}\label{eq:DTA-xim}
\DTA^*(x_{i, m}) = x_{i, m+1}, \hs i \in [1, n], \, m \in \ZZ.
\end{equation}

Turning to $\calX$, recall that $\fX(\calX)$ denotes the set of $y$-variables of $\calX$ and  $\calU(\calX)$
the algebra of regular functions on $\calX$. It is easy to see in examples that $y_1, \ldots, y_n$ are in general not in $\calU(\calX)$.  Recall from \eqref{eq:global-var} that elements in 
$\fX^{\rm global}(\calX) = \fX(\calX) \cap \calU(\calX)$
and called {\it global $y$-variables} of $\calX$. 
Define inductively
\[
Y_i = y_i \prod_{j=1}^{i-1} (1+Y_j)^{-a_{j,i}} \quad \text{and} \quad Y_{n+i} = y_{n+i} \prod_{j=1}^{i-1} (1+Y_{n+j})^{-a_{j,i}}, \hs i \in [1, n].
\]
It is shown in \cite{CdsgL:ringel-conj} that $Y_i$ and $Y_{n+i}$ are global $y$-variables. 
As ${\bf y} =(y_1, \ldots, y_n)$ and ${\bf Y} =(Y_1, \ldots, Y_n)$ are related by invertible positive rational maps, 
the $n$-tuple ${\bf Y}$ gives a coordinate chart on $\calX$ and constitutes global coordinates on $\calXpos$. 
On the other hand, we emphasize that ${\bf Y}$ is in general {\it not} a cluster of $\calX$
(see \cite[Remark 6.4]{CdsgL:ringel-conj} for more details).
It follows from $\DTX^*({\bf y}) = {\bf y}[1]$ that 
\begin{equation}\label{eq:DT-Yi}
\DTX^*(Y_i) = Y_{n+i}, \hs i \in [1,n].
\end{equation}
A direct calculation (see \cite[Proposition 6.2]{CdsgL:ringel-conj}) also gives  
\begin{equation}\label{eq:y-Cox-mutation}
        Y_{i} Y_{n+i} = \prod_{j=i+1}^{n} (1 + Y_j)^{-a_{j,i}} \, \prod_{j=1}^{i-1} (1 + Y_{n+j})^{-a_{j,i}}, \quad i \in [1,n].
\end{equation}
Similar to the case for $\calA$, one has the {\it generic $Y$-frieze pattern} 
which is a surjective map 
$[1, n] \times \ZZ \rightarrow \fX^{\rm global}(\calX),  (i, m) \mapsto Y_{i,m}$ 
such that $Y_{i, 0} = Y_i$ for 
$i \in [1, n]$ and satisfying
\begin{equation}\label{eq:Yim-0}
Y_{i, m}Y_{i,m+1} = \prod_{j = i+1}^n (1+Y_{j,m})^{-a_{j,i}} \, \prod_{j=1}^{i-1} (1+Y_{j,m+1})^{-a_{j,i}}
\end{equation}
for all 
$i \in [1, n]$ and $m \in \ZZ$. See \cite{CdsgL:ringel-conj,Antoine:Y-frieze-patterns}.
An inductive argument using \eqref{eq:DT-Yi} then gives
\begin{equation}\label{eq:DT-Yim}
\DTX(Y_{i,m}) = Y_{i,m+1}, \hs \hs  i \in [1, n], \; m \in \ZZ.
\end{equation}
By \eqref{eq:p-yt}, the ensemble map $p: \calA \to \calX$ satisfies $p^* {\bf y} = {\bf x}^{B_A}$.  
It then follows that
\begin{equation}\label{eq:p-Yi}
p^* Y_i = x_ix_{n+i}-1 = \prod_{j = i+1}^n x_j^{-a_{j,i}} \prod_{j=1}^{i-1} x_{n+j}^{-a_{j,i}}, \hs i \in [1, n].
\end{equation}
By induction, for all $i \in [1,n]$ and $m \in \ZZ$, one has
\begin{equation}\label{eq:p-Yim}
p^* Y_{i,m} = x_{i, m}x_{i, m+1}-1 = \prod_{j=i+1}^n x_{j,m}^{-a_{j,i}} \, \prod_{j=1}^{i-1} x_{j,m+1}^{-a_{j,i}}.
\end{equation}
By \cite[Remark 6.4]{CdsgL:ringel-conj}, the $n$-tuple 
$(x_{1, m+1}, \ldots, x_{i-1, m+1}, x_{i, m}, x_{i+1, m}, \ldots, x_{n, m})$, for 
any $i \in [1, n]$ and  $m \in \ZZ$, is a cluster of $\calA$. Thus 
$p^*Y_{i, m}$, for any $i \in [1, n]$ and  $m \in \ZZ$, is a cluster monomial of $\calA$.

\subsection{Fixed points of cluster DT transformations}\label{ss:fixedpoints} We continue with the
notation from $\S$\ref{ss:finitetype}.  By \leref{:DT-periodic}, $\DTA$ and $\DTX$ are periodic.
 By \cite[Theorem B]{Ishibashi-fxpoint}, $\DTApos:\calApos\to \calApos$ and $\DTXpos:\calXpos\to\calXpos$ admit (potentially multiple) fixed points. Our next result says that such fixed points are unique.

\bth{fixedpoint}
Both $\DTApos: \calApos \to \calApos$ and $\DTXpos: \calXpos\to \calXpos$ have unique fixed points.
\eth

\thref{fixedpoint} follows from the following \prref{:ab-equations} and \leref{:fixpt-lemma2}.

\bpr{:ab-equations}
1) Fixed points of $\DTApos$ in $\calApos$ are in bijection with $(\RR_{>0})^n$-valued solutions to
\begin{equation}
        a_i^2 = 1 + \prod_{j\neq i}a_j^{-a_{j,i}}, \quad \hs  i \in [1,n].\label{eq:fixpt-eq-A}
    \end{equation}

2) Fixed points of $\DTXpos$ in $\calXpos$ are in bijection with $(\RR_{>0})^n$-valued solutions to
\begin{equation}
        b_i^2 = \prod_{j\neq i}(1 + b_j)^{-a_{j,i}}, \quad \hs i \in [1,n].\label{eq:fixpt-eq-X}
    \end{equation}

3) If $(a_1,\ldots, a_n)$ is a solution to \eqref{eq:fixpt-eq-A} then 
$(b_1, \ldots, b_n) =  (a_1^2-1, \, \ldots, \, a_n^2-1)$
is a solution to \eqref{eq:fixpt-eq-X}. 
\epr

\begin{proof}
    Suppose that ${\bf a} \in \calApos$ is a fixed point of $\DTApos$. Then $x_i({\bf a}) = x_{n+i}({\bf a})$ for all $i \in [1,n]$. By \eqref{eq:x-Cox-mutation}, $(a_1, \ldots, a_n) = (x_i({\bf a}), \ldots, x_n({\bf a}))$ is an $(\RR_{>0})^n$-valued solution to the system of equations in \eqref{eq:fixpt-eq-A}. As every point in $\calApos$ is determined by 
    its coordinates in ${\bf x}$, the claim holds for $\calApos$. 

     Suppose that ${\bf b} \in \calXpos$ is a fixed point of $\DTXpos$. Then $Y_i({\bf b}) = Y_{n+i}({\bf b})$ for all $i \in [1,n]$. By \eqref{eq:x-Cox-mutation}, $(b_1, \ldots, b_n) = (Y_1({\bf b}), \ldots, Y_n({\bf b})$ 
     is an $(\RR_{>0})^n$-valued solution to the system of equations \eqref{eq:fixpt-eq-X}. 
     As every point in $\calXpos$ is determined by 
    its coordinates in ${\bf Y}$, the claim holds for $\calXpos$. 
3) is straightforward. 
\end{proof}

\ble{:fixpt-lemma2}
The systems of equations \eqref{eq:fixpt-eq-A} and \eqref{eq:fixpt-eq-X} admit unique $(\RR_{>0})^n$-valued solutions.
\ele
\begin{proof}
    This is well-known in the literature (for \eqref{eq:fixpt-eq-A}, see 
\cite[Theorem 5.3.2]{Lee:fixed-point-Y-system} for the simply-laced case and \cite[Theorem 5.3.4]{Lee:fixed-point-Y-system} for the non simply-laced case via folding; for \eqref{eq:fixpt-eq-X} see \cite[\S 5]{KNS:T-sys}).
We have solved the equations explicitly case-by-case, and the results are given in 
Table \ref{table:A} in Appendix \ref{s:appendixA}. 
\end{proof}

\bre{:fixpt-and-frieze}
By \prref{:ab-equations}, fixed points of $\DTApos$ (resp. of $\DTXpos$) are in one-to-one correspondence with translation-invariant $\RR_{>0}$-valued frieze patterns (resp. $Y$-frieze patterns) associated to $A$ (as defined in \cite{Antoine:Y-frieze-patterns}). An elegant construction of the translation-invariant $\RR_{>0}$-valued frieze pattern of type $A_n$ is given by Conway and Coxeter in \cite{Con-Cox1}. 
\hfill $\diamond$
\ere

We conclude this section with a brief discussion about the relation between $\DTApos$ and the {\it superunitary region}. Recall from \cite{GM:finite} that the latter is the topological subspace $\calA_{\geq 1}$ of $\calApos$ defined by 
\[
   \calA_{\geq 1} = \{a \in \calApos: x(a) \geq 1, \forall x \in \mathfrak{X}(\calA)\},
\]
and its interior is
$\calA_{> 1} = \{a \in \calApos: x(a) > 1, \forall x \in \fX(\calA)\}$.
For $x \in \fX(\calA)$, let $V(x) =  \{ a \in \calApos : x(a) = 1\}$, and set $V = \bigcup_{x \in \fX(\calA)} V(x) \subset \calApos$. Then $\DTApos$ fixes $V$, so $\DTApos$ also fixes $V^\circ:= \calApos \backslash V$. 

\bpr{}
The interior $\calA_{> 1}$ of $\calA_{\geq 1}$ is the only connected component of $V^\circ$ fixed by $\DTApos$,
and the  unique fixed point ${\bf a} \in \calApos$ of $\DTApos$ (see \thref{fixedpoint}) lies in $\calA_{>1}$.
\epr
\begin{proof}
   Since $V^\circ$ has finitely many connected components, it is clear that $\calA_{>1}$ is a connected component fixed by $\DTApos$. 
   Suppose $V_1 \neq \calA_{>1}$ is a connected component of $V^\circ$ fixed by $\DTApos$. Then there exists $x \in \fX(\calA)$ such that $x(p) \in (0,1)$ and $\DTA^*(x)(p) \in (0,1)$ for all $p \in V_1$. Since the map in \eqref{eq:gen-frieze-pat} is surjective, we can find $i \in [1,n]$ and $m \in \ZZ$ such that $x = x_{i,m}$; in particular $\DTA^*(x) = x_{i,m+1}$ by 
   \eqref{eq:DTA-xim}. But then the left-hand side of \eqref{eq:frieze-mut-0} is strictly smaller than 1 whilst the right-hand side of \eqref{eq:frieze-mut-0} is strictly greater than 1, which is a contradiction. 
   To see that ${\bf a} \in \calA_{>1}$, 
note that by \prref{:ab-equations}, $x_{i, 1}({\bf a}) = x_i({\bf a}) >1$ for every $i \in [1, n]$, so by
\eqref{eq:DTA-xim}, $x({\bf a}) > 1$ for every $x \in \fX(\calA)$. Thus ${\bf a} \in \calA_{>1}$.
\end{proof}

\bre{:fixpt-superunitary}
It was shown in \cite[Theorem A]{GM:finite} that $\calA_{\geq 1}$ is homeomorphic to the closed ball $\mathbb{D}_n$. Since $\DTApos$ restricts to a homeomorphism of $\calA_{\geq 1}$
to itself, the Brouwer fixed point theorem \cite[p.65]{GuilleminPollack:Diff-Top} implies the {\it existence} of a fixed point in $\calA_{\geq 1}$ (which by the arguments above, lies in $\calA_{>1}$). This yields an alternative proof of the existence part of 1) in \thref{fixedpoint}.  
\ere

\subsection{Cluster exponents}\label{ss:cluster exponents}
Continuing with the notation of \S \ref{ss:fixedpoints}, let ${\bf a} \in \calApos$ and ${\bf b} \in \calXpos$
be the respective fixed points of  $\DTApos$ and $\DTXpos$,  and consider the linear maps 
\[
\dtA ={\rm d}_{\bf a}\DTApos: \;{\rm T}_{\bf a} (\calApos) \longrightarrow {\rm T}_{\bf a} (\calApos)
\hs \mbox{and} \hs 
\dtX ={\rm d}_{\bf b}\DTXpos:\; {\rm T}_{\bf b} (\calXpos) \longrightarrow {\rm T}_{\bf b} (\calXpos).
\]
By \thref{:same-char}, $\dtA$ and $\dtX$ have the same characteristic polynomial which we denote by $\calP_\Phi(x)$. 

\bre{:fpIsLefschetz}
{\rm
The points ${\bf a} \in \calApos$ and ${\bf b}\in \calXpos$ are Lefschetz  fixed points of $\DTApos$ and $\DTXpos$ respectively, in the sense 
\cite[p.33]{GuilleminPollack:Diff-Top}
that $\dtA$ and $\dtX$ do not have eigenvalue 1. 
Indeed, the action on $\calApos$ of the finite subgroup $H$ of Diff($\calApos$) generated by $\DTApos$ is
linearizable at ${\bf a}$ (see, for example, \cite[Prop. 5, no.3, \S 9]{bourbaki:LieAlgLieGp-9}), and an eigenvector 
of $\dtA$ with eigenvalue $1$ would give finitely many fixed points of $\DTApos$ in $\calApos$ near ${\bf a}$.
Similarly, $\dtX$ does not have eigenvalue $1$.
\hfill $\diamond$
}
\ere

Let again ${\bf h}$ be the Coxeter number of $\Phi$. 
By \deref{:linear-expo-intro} and \deref{:cluster-expo-intro}, the exponents of the linear maps 
$\dtA$ and $\dtX$ with respect to ${\bf h} + 2$ are called  the cluster exponents of the cluster ensemble $(\calA,\calX, p)$ with respect to ${\bf h}+2$.

\bre{:relation-with-Y-sys-exponents}
When $\Phi$ is simply laced (i.e. of type $ADE$), the exponents of $\dtX$ are precisely the {\it exponents of the quiver $Q(\Phi,\ell = 2)$}, defined by Y. Mizuno (see \cite[\S 3.5]{Mizuno:Y-sys-exponents} for the notation). We emphasize that when $\Phi$ is {\it not} simply laced, the exponents of $\dtX$ and the exponents of $Q(\Phi, 2)$ are different, owing to the fact that the cluster ensembles underlying either one are distinct.
\ere

To compute the cluster exponents of $(\calA, \calX, p)$, we first give matrix representations 
of $\dtA^*$ and $\dtX^*$ in suitable bases of the co-tangent spaces $T_{\bf a}^*(\calApos)$ and 
$T_{\bf b}^*(\calXpos)$. Set again $a_i = x_i({\bf a})$, $b_i = Y_i({\bf b})$, and 
\[
    v_i = \frac{d_{\bf a} x_i}{a_i} \in T_{\bf a}^* (\calApos) \quad \text{and} \quad w_i = \frac{d_{\bf b}Y_i}{1+b_i}
    \in T_{\bf b}^* (\calXpos), \hs i \in [1, n].
\]
Then  
$(v_1,\ldots, v_n)$ is a basis of $T_{\bf a}^*\calApos$ and $(w_1, \ldots , w_n)$ is a basis of $T_{\bf b}^*\calXpos$.
Set
\begin{equation}\label{eq:ci}
    \kappa_i = \frac{a_i^2}{a_i^2 -1} = \prod_{j=1}^n a_j^{a_{j,i}} = \frac{b_i+1}{b_i},
    \hs i \in [1, n]
\end{equation}
(recall \eqref{eq:fixpt-eq-A} and  3) of \prref{:ab-equations}). Write $\kappa = (\kappa_1, \ldots, \kappa_n) \in (\RR_{>1})^n$, so  $\kappa = (a_1,\ldots, a_n)^A$. 
Recall $L_\kappa$ and $U_\kappa$ defined in \eqref{eq:LU-intro}. Set again $M_\kappa = -L_\kappa U_\kappa^{-1}$
as in \eqref{eq:Mc}.

\bpr{:matrix-realization-linearDT}
With the notation as above, we have 
\begin{equation}
    \dtA^* (v_1,\ldots, v_n) = (v_1,\ldots, v_n) M_{\kappa} \quad \text{ and } \quad \dtX^* (w_1,\ldots, w_n) = (w_1,\ldots, w_n) M_{\kappa}.
\end{equation}
Consequently, the common characteristic polynomial $\calP_\Phi(x) \in \RR[x]$ of 
 $\dtA$ and 
$\dtX$ is given by 
\[
\calP_\Phi(x) = \det(x-M_\kappa) = \frac{1}{\kappa_1\kappa_2\cdots \kappa_n} \det (xU_{\kappa}+L_{\kappa}).
\]
\epr

\begin{proof}
Recall that $ (x_{n+1}, \ldots, x_{2n}) = \DTA^*(x_1,\ldots, x_n)$. 
For $i \in [1,n]$, we then have
$\dtA^*(d_{\bf a}x_i) = d_{\bf a}x_{n+i}$. 
Applying the differential $d$ to both sides of \eqref{eq:x-Cox-mutation} for each $i \in [1,n]$
and using \eqref{eq:d-xV-1}, we get
    \[
x_ix_{n+i}\left(\frac{dx_i}{x_i} + \frac{dx_{n+i}}{x_{n+i}}\right) = \left(\prod_{j=i+1}^n x_j^{-a_{j, i}} \prod_{j=1}^{i-1} x_{n+j}^{-a_{j, i}}\right) \left(
\Sigma_{j=i+1}^n \left(-a_{j, i} \frac{dx_j}{x_j}\right) + \Sigma_{j=1}^{i-1} \left(-a_{j, i} \frac{dx_{n+j}}{x_{n+j}}\right)\right).
\]
Evaluating at  ${\bf a} \in \calApos$ and using the fact that $a_i^2 - 1 = \prod_{j\neq i}a_j^{-a_{j,i}}$, we have
\begin{equation}\label{eq:delta-1}
\frac{a_i^2}{a_i^2-1}(v_i+\dtA^*(v_i)) = \Sigma_{j=i+1}^n \left(-a_{j, i} v_j\right) + 
\Sigma_{j=1}^{i-1} \left(-a_{j, i} \dtA^*(v_j)\right), \hs i \in [1, n]. 
\end{equation}
Using \eqref{eq:ci} and \eqref{eq:LU-intro} to rewriting \eqref{eq:delta-1}, we have $(\dtA^*(v_1), \ldots, \dtA^*(v_n))U_{\kappa} = -(v_1, \ldots, v_n)L_{\kappa}$. Thus
\begin{equation}\label{eq:delta-u}
\dtA^*(v_1, \ldots, v_n) = (v_1, \ldots, v_n) M_{\kappa}.
\end{equation}

Turning to $\dtX^*$, recall that  $ (Y_{n+1}, \ldots, Y_{2n}) = \DTX^*(Y_1,\ldots, Y_n)$, so
$\dtX^*\, (d_{\bf b}Y_i) = d_{\bf b}Y_{n+i}$ for $i \in [1,n]$.
Applying the differential $d$ to both sides of \eqref{eq:y-Cox-mutation} for each $i \in [1,n]$ and evaluating at ${\bf b}$, 
we get
    \begin{equation}\label{eq:Y-twist-linear}
       (1+b_i)b_i (w_i + \dtX^*(w_i)) = \prod_{j\neq i} (1+b_j)^{-a_{j,i}} \left( \sum_{j=i+1}^n (-a_{j,i})w_j + \sum_{j=1}^{i-1}(-a_{j,i})\dtX^*(w_j)\right).
    \end{equation}
Using \eqref{eq:ci} and  \eqref{eq:fixpt-eq-X} to rewrite \eqref{eq:Y-twist-linear}, we get 
    \[
        \kappa_iw_i + \sum_{j=i+1}^n a_{j,i} w_j = - \left( 
        \kappa_i\dtX^*(w_i) + \sum_{j=1}^{i-1} a_{j,i} \dtX^*(w_j)\right), 
    \]
    and thus $\dtX^*(w_1, \ldots, w_n) = (w_1, \ldots, w_n) M_{\kappa}$, 
    which concludes the proof.
\end{proof}

\bre{:u-variables}
\prref{:matrix-realization-linearDT} highlights the importance of the vector 
$\kappa =(\kappa_1,\ldots, \kappa_n) \in (\RR_{>1})^n$. 

i) Setting $u_{i, m} = \frac{Y_{i,m}}{1+Y_{i, m}}$ for $i \in [1, n]$ and $m \in \ZZ$,
so that $Y_{i, m} = \frac{u_{i,m}}{1-u_{i,m}}$, the relations in \eqref{eq:Yim-0} become
\begin{equation}\label{eq:local-u-eq}
    \frac{u_{i,m}}{1-u_{i,m}} \, \frac{u_{i,m+1}}{1-u_{i,m+1}} = \prod_{j = i+1}^n (1-u_{j,m})^{a_{j,i}} \, \prod_{j=1}^{i-1} (1-u_{j,m+1})^{a_{j,i}}.
\end{equation}
By \eqref{eq:DT-Yim}, $\DTX^*(u_{i,m}) = u_{i,m+1}$ for all $i,m$. 
The $u_{i,m}$'s thus defined are called {\it u-variables} and the relations \eqref{eq:local-u-eq} the {\it local u-equations} in \cite[$\S$10.2]{Arkani-hamed:cluster-config} and \cite[p.6]{Arkani-Hamed-binary}. The
tuple $(u_i = u_{i, 0}({\bf b}))_{i \in [1,n]} = (\kappa_i^{-1})_{i \in [1,n]}$ is the unique solution  in $(0,1)^n$ to the system \cite{Nahm:CFT,Vlasenko:NahmConjecture} and \cite[Proposition 4.2.3]{Lee:fixed-point-Y-system}
\begin{equation}\label{eq:u-system}
    u_i^2 = \prod_{j=1}^n (1-u_j)^{a_{j,i}}, \quad i \in [1,n].
\end{equation}

ii) Set $\log(\kappa) = (\log(\kappa_1),\ldots, \log(\kappa_n)) >0$ and 
$\log(a) = (\log(a_1),\ldots, \log(a_n))>0$. Note then that 
\begin{equation}\label{eq:vinberg-crit}
    \log(\kappa) = \log(a)A.
\end{equation}
Recall that an integer $n \times n$ matrix $A = (a_{i,j})$ is called a {\it symmetrizable generalized 
Cartan matrix} if 
(1) $a_{i,i} = 2$ for all $i \in [1,n]$, (2) $a_{i,j} \leq 0$ for all $i \neq j$, and (3) there exists a diagonal matrix $D$ with positive integer diagonal entries such that $DA$ is symmetric. 
By Vinberg's criterion \cite[Corollary 4.3]{Kac:inf-dim-Lie-alg}, if $A$ is an $n \times n$ 
symmetrizable generalized 
Cartan matrix and if there exists $v \in (\RR_{>0})^n$ such that $Av \in (\RR_{>0})^n$, then $A$ must be of 
finite type. We can thus conclude that if the DT transformations of an {\it acyclic cluster 
ensemble} $(\calA,\calX,p)$ admit positive fixed points, then the cluster ensemble must be of finite type. 
Indeed, the matrix pattern of such an $(\calA,\calX,p)$  contains a matrix which is, up to simultaneous re-labeling or rows and columns, is of the form $B_A$ for a  
symmetrizable generalized Cartan matrix $A$. The same arguments presented in \S \ref{ss:finitetype} and \prref{:ab-equations} show that $(\calA,\calX,p)$ admits DT transformations, and fixed points of $\DTApos$ and $\DTXpos$ 
are in one-to-one correspondence with $\mathbb{R}_{>0}$-valued solutions to \eqref{eq:fixpt-eq-A} and \eqref{eq:fixpt-eq-X} respectively. If $a = (a_1,\ldots, a_n)$ is a solution to \eqref{eq:fixpt-eq-A}, then defining $\kappa$ as in \eqref{eq:ci}, we obtain that $\kappa$ and $a$ satisfy \eqref{eq:vinberg-crit}, and thus $A$ is of finite type by Vinberg's criterion.
\hfill $\diamond$
\ere

Turning now to the explicit formula for the polynomial $\calP_\Phi(x)$, consider the Weyl group $W$ of $\Phi$,
and let 
$d_1  \leq \cdots \leq  d_n$
be {\it degrees of fundamental invariants} of $W$, or simply the {\it degrees of $W$}.
For the convenience of the reader, we recall in 
Table \ref{table:d} in Appendix B the explicit values of the degrees in all types from 
\cite[Page 59]{HumphreysCoxGrp}. Using the fact \cite[Theorem 3.19]{HumphreysCoxGrp} that
$d_1-1, \ldots, d_n-1$
are the exponents of any Coxeter element of $W$ with respect to the Coxeter number ${\bf h}$, one has
$2 = d_1 \leq d_2 \leq \cdots \leq d_n = {\bf h}$ and 
\begin{equation}\label{eq:di-di}
d_i + d_{n+1-i} = {\bf h} + 2.
\end{equation}
Set $\zeta = e^{\frac{2\pi \sqrt{-1}}{{\bf h}+2}}$ and define
\[
    \calD_\Phi(x) = \prod_{j=1}^n (x - \zeta^{d_j}).
\]

\bth{:expnts-conjecture} For any indecomposable finite root system $\Phi$
and any cluster ensemble $(\calA, \calX, p)$ 
of type $\Phi$, one has 
$\calP_{\Phi}(x) = \calD_\Phi(x)$, 
where $\calP_{\Phi}(x)$ is the (common) characteristic polynomial of $\dtA$ and $\dtX$.
\eth

The proof of \thref{:expnts-conjecture}, which will be given in Appendix \ref{s:appendixB}, 
is achieved in a case-by-case analysis, by explicitly computing 
$\calP_\Phi(x)$ using  \prref{:matrix-realization-linearDT}.

\subsection{A polynomial identity for root systems and Mizuno's conjecture}\label{ss:applications}
Continuing with the notation in $\S$\ref{ss:cluster exponents}, 
choose any system of simple roots  $\{\alpha_1,\ldots, \alpha_n\}$ of $\Phi$ and define the
{\it height function} 
\[
 {\rm ht}: \;\;\Phi \longrightarrow \ZZ,\;\;   {\rm ht}\left( \sum_{i = 1}^n h_i \alpha_i \right ) = \sum_{i = 1}^n h_i.
\]

\bpr{:degrees-and-heights}
With $\zeta = e^{\frac{2\pi \sqrt{-1}}{{\bf h}+2}}$, one has\footnote{We expect the identity in 
\eqref{eq:degrees-and-heights} to be known, but we are unable to find a reference in the literature.} 
    \begin{equation}\label{eq:degrees-and-heights}
        \left(\frac{x^{{\bf h}+2}-1}{x-1}\right)^n = \calD_\Phi(x)  \prod_{\beta \in \Phi} \left(x- \zeta^{{\rm ht}(\beta)}\right).
    \end{equation}
\epr

\begin{proof} Let $R = [1, n] \times [1, {\bf h}+1]$, and for $j \in [1, {\bf h}+1]$, let
$P_{1, j}(x)  = \cdots = P_{n, j}(x) =x - \zeta^j$.
Then
\begin{align*}
\left(\frac{x^{{\bf h}+2}-1}{x-1}\right)^n &= \prod_{j=1}^{{\bf h}+1} (x-\zeta^j)^n = 
\prod_{(i, j) \in R} P_{i, j}(x) =\prod_{(i, j) \in R_1} P_{i, j}(x) \prod_{(i, j) \in R_2} P_{i, j}(x) 
\prod_{(i, j) \in R_3} P_{i, j}(x),
\end{align*}
where  $R_1 = \{(i, j) \in R: d_i = j\}$, $R_2 = \{(i, j) \in R: d_i \geq  j+1\}$, and
$R_3 = \{(i, j) \in R: d_i \leq  j-1\}$.
Note that 
\[
\prod_{(i, j) \in R_1} P_{i, j}(x) = \prod_{i \in [1, n]} (x-\zeta^{d_i}) = \calD_\Phi(x).
\]
Let $\Phi_+ \subset \Phi$ be the set of positive roots in $\Phi$,
and recall that ${\rm ht}(\beta) \in [1, {\bf h}-1]$ for every $\beta \in \Phi_+$. 
By a theorem of Kostant (see \cite[Theorem 3.20]{HumphreysCoxGrp}), for each $j \in [1, {\bf h}-1]$ one has
\begin{equation}\label{eq:fj}
|\{i \in [1, n]: d_i \geq j+1\}|  = 
|\{\beta \in \Phi_+:\; {\rm ht}(\beta) = j\}|.
\end{equation}
Note that \eqref{eq:fj} also holds for $j = {\bf h}$ or $j \in {\bf h} +1$ as both sides are 
$0$ in these cases. Thus
\[
|\{i \in [1, n]: (i, j) \in R_2\}| = |\{\beta \in \Phi_+:\; {\rm ht}(\beta) = j\}|
\]
for every $j \in [1, {\bf h}+1]$. It thus follows that 
\[
\prod_{(i, j) \in R_2} P_{i, j}(x) = 
\prod_{j \in [1, {\bf h}+1]} (x-\zeta^j)^{|\{i \in [1, n]: \;(i, j) \in R_2\}|}=
\prod_{j \in [1, {\bf h}+1]} (x-\zeta^j)^{|\{\beta \in \Phi_+: \;{\rm ht}(\beta) = j\}|}
=\prod_{\beta \in \Phi_+}(x-\zeta^{{\rm ht}(\beta)}).
\]
By \eqref{eq:di-di}, for $(i, j) \in R$, one has 
\[
(i, j) \in R_3 \; \Longleftrightarrow \; d_i \leq j-1\; \Longleftrightarrow \;
d_{n+1-i} \geq {\bf h}+2 - (j -1) \; \Longleftrightarrow \; (n+1-i, \; {\bf h}+2-j) \in R_2.
\]
It follows that
\[
\prod_{(i, j) \in R_3} P_{i, j}(x) = \prod_{(n+1-i, \; {\bf h}+2-j) \in R_2}(x-\zeta^j)
= \prod_{(i, j) \in R_2} (x-\zeta^{{\bf h}+2-j})=\prod_{(i, j) \in R_2}(x-\zeta^{-j})=
\prod_{\beta \in -\Phi_+}(x-\zeta^{{\rm ht}(\beta)}).
\]
This proves \eqref{eq:degrees-and-heights}.
\end{proof}

\bco{:Mizuno-conj}
Mizuno's conjecture \cite[Conjecture 3.8]{Mizuno:Y-sys-exponents} holds for $\Phi$ of type $ADE$ and $\ell = 2$. 
\eco

\begin{proof}
Mizuno's conjecture \cite[Conjecture 3.8]{Mizuno:Y-sys-exponents} for $\Phi$ of simply-laced 
type and $\ell = 2$ 
says precisely that \eqref{eq:degrees-and-heights} holds when $\calD_\Phi(x)$ is replaced by $\calP_\Phi(x)$. 
Since $\calP_\Phi(x) = \calD_\Phi(x)$ by \thref{:expnts-conjecture}, the conjecture follows 
from \prref{:degrees-and-heights}.
\end{proof}


\bre{:mizuno-conjecture-proofs}
Mizuno's conjecture \cite[Conjecture 3.8]{Mizuno:Y-sys-exponents} for type $A_n$ was verified by himself
in \cite{Mizuno:Y-sys-exponents}, and the type $D_n$ case  was recently verified by R. Takenaka \cite{Takenaka-Dn}.  
\ere

\section{An affine model via double Bruhat cells}\label{s:Lie}
\subsection{Lie theory preliminaries}\label{ss:rdBc-and-twist} Let $\Phi$ be an irreducible finite root system. 
Let $G$ be a connected and simply connected simple Lie group over $\CC$ of type $\Phi$, 
let $(B,B_-)$ be a pair of opposite Borel subgroups of $G$, and let $T = B \cap B_-$,
a maximal torus of $G$. Let $\mathfrak{t}$ be the Lie algebra of $T$, and identify $\Phi$ with the set of roots in $\mathfrak{t}^*$ defined by $(G,T)$. Let 
$\{\alpha_1,\ldots, \alpha_n\}\subset \Phi_+ \subset \Phi$ be the respective sets of simple roots and positive roots determined by $B$, let $\{\alpha_1^\vee, \ldots, \alpha_n^\vee \}\subset \t$ be the corresponding  {\it simple co-roots}, and let 
$A = (a_{i, j})$ where 
$a_{i,j} = \alpha_i^\vee(\alpha_j)$ for $i, j \in [1, n]$.
Let $\{\omega_1, \ldots, \omega_n\}$ be the set of {\it fundamental weights}, i.e., 
the basis of $\t^*$ that is dual to the basis
$\{\alpha_1^\vee, \ldots, \alpha_n^\vee\}$ of $\t$. Writing $\underline{\alpha} = (\alpha_1, \ldots, \alpha_n)$
and $\underline{\omega} = (\omega_1, \ldots, \omega_n)$, we have
\[
\underline{\alpha} = \underline{\omega} A.
\]
Let $\calP =\sum_{i=1}^n \ZZ \omega_i \subset \t^*$ be the weight lattice. 
Define $\mu > \delta$ for $\mu,\delta \in \calP$ if $\mu- \delta \in \ZZ_{\geq 0}\alpha_1 + \cdots + \ZZ_{\geq 0}\alpha_n$.

Let $W$ 
be the Weyl group of $(G,T)$ with generating set $S = \{s_1, \ldots, s_n\}$, where $s_i \in W$ 
for $i \in [1, n]$
is the {\it simple reflection} associated to $\alpha_i$. For $i \in [1,n]$, we fix root vectors $e_i$ of $\alpha_i$ and $e_{-i}$ of $-\alpha_i$ such that $[e_i,e_{-i}] = \alpha_i^\vee$. The triple $\{e_i,e_{-i},\alpha_i^\vee\}$ determines a Lie group embedding $\theta_i: {\rm SL}_2 \rightarrow G$. The element 
$\begin{pmatrix}    0 & -1 \\ 1 & 0
    \end{pmatrix}$
is a representative of $s_i$ in the normalizer 
$N_G(T)$ of $T$ in $G$. Each $w \in W$ then has the well-defined representative $\overline{w}$  in $N_G(T)$, determined by  the relation $\overline{u_1 \, u_2} = \overline{u}_1 \, \overline{u}_2$ whenever $\ell (u_1 u_2) = \ell (u_1) + \ell (u_2)$, where $\ell: W \to \NN$ is the 
length function with respect to the Coxeter system $(W,S)$. Let $w_0 \in W$ be the unique element of maximal length in $W$. For $t \in T$ and $w \in W$, set
$t^w = \overline{w}^{\, -1} t \overline{w} \in T$.

Let $N$ and $N_-$ be the respective unipotent radicals of $B$ and $B_-$. Then every $g \in B_-B$ can be
uniquely decomposed as $g = [g]_-\, [g]_0\,[g]_+$ for unique $([g]_-, [g]_0,[g]_+) \in N_- \times T \times N$.
By \cite{FZ:double}, for  $u, v \in W$ and $i \in [1,n]$, one has the regular function 
$\Delta_{u\omega_i, \,v\omega_i}$ on $G$, called a {\it generalized minor},  
such that 
\[
\Delta_{\omega_i, \omega_i}(g) = [g]_0^{\omega_i} \; \;\;\mbox{for}\; \g \in B_-B\hs 
\mbox{and} \hs \Delta_{u\omega_i,\,v\omega_i}(g) =\Delta_{\omega_i, \, \omega_i} (\overline{u}^{\, -1} g \overline{v}) \;\;\; \mbox{for}\;\;  g \in G.
\]
For $i \in [1, n]$, let
$i^* \in [1, n]$ be such that $w_0\alpha_i = -\alpha_{i^*}$, where $w_0$  is the longest element of $W$.
By  \cite[(2.15) and (2.25)]{FZ:double},
for every $i \in [1, r]$  and for all $g \in G$, one has
\begin{equation}\label{eq:w0-g}
\Delta_{w_0\omega_{i^*}, w_0\omega_{i^*}}(g^{-1}) = \Delta_{\omega_i, \omega_i}(g).
\end{equation}

We will consider the Coxeter element $c = s_1s_2 \cdots s_n \in W$. 
For $i \in [1, n]$, let $\beta_i = s_1\cdots s_{i-1}\alpha_i$, and write
$\underline{\beta} = (\beta_1, \ldots, \beta_n)$. It follows from the definition that 
$\underline{\beta} = \underline{\omega}-c\,\underline{\omega}$ on the one hand, and $\underline{\beta} 
U_{\bf 1}= 
\underline{\alpha}$, on the other, where $U_{\bf 1}$ is given in \eqref{eq:LU-intro}
for ${\bf 1} = (1, \ldots, 1)$. With $M_{\bf 1} = -L_{\bf 1} U_{\bf 1}^{-1}$ as in \eqref{eq:Mc}, one thus has
\begin{equation}\label{eq:M-one}
c\,\underline{\omega} = \underline{\omega} M_{\bf 1}.
\end{equation}

\bre{:c}
{\rm
Note that $(1-c)\, \underline{\omega} =\underline{\beta} = \underline{\alpha} U_{\bf 1}^{-1}$, from which we 
can deduce  (the well-known facts)
that $c$, as a linear map on $\t^*$, does not have $1$ as an eigenvalue, and that $1-c: \calP \to \calQ$ is a lattice isomorphism, where $\calQ = \sum_{i=1}^n \ZZ \alpha_i$ is the root lattice.
\hfill $\diamond$
}
\ere

\subsection{The Yang-Zelevinsky cluster structure on $\Lcc$}\label{ss:Lcc} 
With 
$c = s_1s_2\cdots s_n \in W$, consider now the  
{\it reduced double Bruhat cell}
$\Lcc = N \overline{c} N \cap B_-c^{-1}B_- \subset G$.
By \cite[Theorem 1.1]{FZ:double} $\Lcc$ is 
a smooth affine variety of dimension $2n$ (see also \cite[Theorem 7.5]{GHL:friezes}). 

By \cite[Proposition 1.3]{YZ:Lcc},  there is an integer $h(i; c) \geq 1$ for each $i \in [1, n]$ such that
\[
\omega_i > c\omega_i > c^2 \omega_i > \cdots > c^{h(i; c)} \omega_i = w_0 \omega_i = -w_{i^*}.
\]
Following \cite{YZ:Lcc}, let $\Pi(c)=\{c^m\omega_i: \,i \in [1, n], m \in [0, h(i;c)]\}$, and for
$\gamma \in \Pi(c)$ let $u_{\gamma} = \Delta_{\gamma, \gamma}|_{\Lcc}$. Set
\[
f_{i} = \Delta_{\omega_i, c\omega_i}|_{\Lcc} \hs \mbox{and} \hs
p_{i} = f_{i}
\prod_{j=1}^{i-1} f_{j}^{a_{j, i}} \in \CC[\Lcc], \hs i \in [1, n],
\]
and let ${\bf f} = (f_{1}, \ldots, f_{n})$, ${\bf p} = (p_{1}, \ldots, p_{n})={\bf f}^{U_{\bf 1}}$,
and ${\bf u}_{-\underline{\omega}} = (u_{-\omega_1}, \ldots, u_{-\omega_n})$. 
Consider the seed\footnote{The initial seed used in \cite{YZ:Lcc} is 
$((u_{\omega_1}, \ldots, u_{\omega_n}, {\bf p}_c), (-B_A\pp I_n))$ which becomes $\Sigma^{\Lcc}$
after mutations along the 
sequence $(n, \ldots, 2, 1)$.}
\[
\Sigma^{\Lcc} =  \left(({\bf u}_{-\underline{\omega}}, \, {\bf p}_c), \; (B_A\pp I_n)\right)
\]
in $\CC(\Lcc)$ of rank $n$ with ${\bf p}$ as frozen variables.

\bld{:YZ} \cite[Theorem 1.2]{YZ:Lcc}
{\rm The seed $\Sigma^{\Lcc}$ defines a 
cluster structure on 
$\Lcc$, which we call the {\it Yang-Zelevinsky cluster structure} on $\Lcc$,
whose set of cluster variables is $\{u_{\gamma}: \gamma \in \Pi(c)\}$.
}
\eld

Note that the Yang-Zelevinsky cluster structure on $\Lcc$ has principal coefficients at the seed 
$\Sigma^{\Lcc}$. Identifying $\ZZ^n \cong \calP$ using the basis of the fundamental weights $\underline{\omega}
=(\omega_1, \ldots, \omega_n)$, it is shown in \cite[Theorem 1.10]{YZ:Lcc} that 
the $g$-vector of $u_{\gamma}$ with respect to the seed $\Sigma^{\Lcc}$ is $\gamma$.
By \cite[(1.10)]{YZ:Lcc}, 
\begin{equation}\label{eq:u-omega-omega}
u_{-\omega_i} u_{\omega_i} = 1 +p_{i} \prod_{j=i+1}^n u_{-\omega_j}^{-a_{j, i}}
\prod_{j=1}^{i-1} u_{\omega_j}^{-a_{j, i}}, \hs i \in [1, n].
\end{equation}
Setting ${\bf u}_{\underline{\omega}} = (u_{\omega_1}, \ldots, u_{\omega_n})$, it
then follows from \eqref{eq:u-omega-omega} and a direct calculation that 
\begin{equation}\label{eq:mu-n1}
\mu_n \cdots \mu_2\mu_1 \left(\Sigma^{\Lcc}\right) = ({\bf u}_{\underline{\omega}}, \; (B_A\pp -I_n)).
\end{equation}
Introduce the affine sub-variety 
$\Acc =\{g \in \Lcc: \;p_{i;c}(g) = 1, i \in [1,n]\}$ of $\Lcc$, and set  
\[
x_{\gamma;c} = \Delta_{\gamma,\gamma}\vert_{\Acc}, \hs \gamma \in \Pi(c).
\]
 Then
the seed $\Sigma^{\Acc} = ({\bf x}_{-\underline{\omega}} :=(x_{-\omega_i}, \ldots, x_{-\omega_n}), B_A)$ in $\CC(\Acc)$ defines a 
cluster algebra structure on $\Acc$ with trivial 
coefficients whose set of cluster variables is $\{x_{\gamma}: \gamma \in \Pi(c)\}$. 

Recall now the setup of \S \ref{s:fixedpoint} and consider the cluster ensemble $(\calA, \calX, p)$ 
of type $\Phi$ with trivial coefficients. Then $\Acc$ as an affine model of
$\calA$  in the sense that $\Acc$ is affine and that $\CC[\Acc] \cong \calU(\calA) \otimes \CC$. 
The set $\calA(\CC)$ of $\CC$-points of $\calA$ can be identified with 
the union of all the cluster tori in $\Acc$. The complement of $\calA(\CC)$ in 
$\Acc$, called the {\it deep locus} of $\calA$ \cite{deepPt-Muller}, is in general not empty but has co-dimension at least $2$. Similarly, $\Lcc$ is an affine model for $\Aprin$. 
Comparing the seeds $\Sigma^{\Lcc}$ and 
$\Sigma^{\Acc}$, we see that $\Lcc$, with the Yang-Zelevinsky cluster structure, is a principal extension of $\Acc$ at the seed $\Sigma^{\Acc}$. 

For an affine model for $\calX$, consider the action of $T$ on $\Lcc$ given by 
\[
t \cdot g = t g (t^{-1})^c, \hs t \in T, \; g \in \Lcc.
\]
It follows from the definitions that for each $\gamma \in \Pi(c)$ and for all $t \in T$ and $g \in \Lcc$, one has
\[
u_{\gamma} (t\cdot g) = t^{\gamma - c\gamma} u_{\gamma}(g).
\]
Thus the 
$T$-weight of $u_{\gamma}$ is $(1-c)\gamma$ for $\gamma \in \Pi(c)$. 
Let $Z(G)$ be the center of $G$ and let $T_{\rm ad} = T/Z(G)$. Then the character lattice of $T_{\rm ad}$
is the root lattice $\calQ$ which we identify with $\calP$ via the lattice isomorphism 
$1-c: \calP \to \calQ$ (see \reref{:c}), and the action of $T$ on $\Lcc$ descends to an action of $T_{\rm ad}$ on $\Lcc$. 
By Hilbert's finiteness theorem, $\CC[\Lcc]^{T_{\rm ad}}$ is finitely generated, so we have the affine
variety
\[
\Xcc =\Lcc /\!/ T_{\rm ad} := {\rm Spec}(\CC[\Lcc]^{T_{\rm ad}}).
\]
Let $p^{c, c^{-1}}: \Lcc \to \Xcc$ be the projection induced by the algebra embedding
$\CC[\Lcc]^{T_{\rm ad}} \to \CC[\Lcc]$.  By \cite[Lemma 7.10]{GHKK}, we can regard $\Xcc$ as an affine model of $\calX$ in the sense that $\Xcc$ is affine and that $\CC[\Xcc] \cong \calU(\calX) \otimes \CC$. 
Denote also by $p^{c, c^{-1}}: \Acc \to \Xcc$ the restriction of $p^{c, c^{-1}}: \Lcc \to \Xcc$ to $\Acc \subset \Lcc$. 
We then have the affine models
\[
(\Acc, \; \Xcc, \; p^{c, c^{-1}}) \hs \mbox{and} \hs (\Lcc, \; \Xcc, \; p^{c, c^{-1}}),
\]
of the cluster ensembles $(\calA, \calX, p)$ and its principal extension  $(\Aprin, \calX, p^{\rm prin})$, respectively.


\subsection{The Fomin-Zelevinsky twists as DT transformations}\label{ss:FZ-twist}
Following \cite{BZ:tensor}, 
define the {\it Fomin-Zelevinsky twist} on $\Lcc$
to be the biregular isomorphism 
\begin{equation}\label{eq:FZ-twist}
\psi: \;\Lcc \longrightarrow \Lcc,\;\; \,\psi(g)
= [\overline{c}^{\, -1} g]_-^{-1} \overline{c}^{\, -1}  g \, \overline{c} \, [g \overline{c}]_+^{-1}, \quad g \in \Lcc.
\end{equation}
It follows from \eqref{eq:FZ-twist} that $\psi(\Acc) =\Acc$. Set 
${\phi}=\psi|_{\Acc}: \Acc \rightarrow \Acc$.

On the other hand, one checks from the definitions that  $\psi(t \cdot g) = t^c \cdot g$ for all $t \in T$ and $g \in \Lcc$.
Thus $\psi: \Lcc \to \Lcc$ descends to a well-defined bi-regular automorphism
$[\psi]:  \Xcc \rightarrow \Xcc$.

Recall again the matrices $L_{\bf 1}$ and $U_{\bf 1}$ given in \eqref{eq:LU-intro}. 

\bpr{:psi-c} With respect to the Yang-Zelevinsky cluster structure on $\Lcc$, the
Fomin-Zelevinsky twist $\psi: \Lcc \to \Lcc$ is a twisted 
automorphism of $\Lcc$ of DT-type,   and 
\begin{equation}\label{eq:FZtwist-u-p}
\psi^*({\bf u}_{-\underline{\omega}}) = {\bf u}_{\underline{\omega}} \,{\bf p}^{-U_{\bf 1}^{-1}}, \hs \hs 
\psi^*({\bf p}) = {\bf p}^{-U_{\bf 1}^{-1}L_{\bf 1}}.
\end{equation}
Consequently, ${\phi}$ is the DT transformation of $\Acc$,  and $[\psi]$ is the DT transformation of $\Xcc$.
\epr

\begin{proof} Let $N_c = N \cap \oc N_- \oc^{\, -1}$ and $N_{-, c} = N_- \cap \oc N \oc^{\, -1}$. Let $g \in \Lcc$ and write $g$ uniquely as 
\begin{equation}\label{eq:g-decomp}
g = n_c \oc\, n' = m' t\, \oc^{\, -1} m_c,
\end{equation}
where 
$n_c \in N_c, n' \in N, m' \in N_-, t \in T$ and $m_c \in N_{-, c}$. Then
\[
\psi(g) =n'   m_c^{-1} \oc = (n_c\oc)^{-1} m't.
\]
For $i \in [1, n]$, one then has
$f_i(g) = \Delta_{\omega_i, c\omega_i}(g) = \Delta_{\omega_i, \omega_i}(m' t\, \oc^{\, -1} m_c \oc) = t^{\omega_i}$,
and it  follows from 
\eqref{eq:w0-g} that
\begin{align*}
(\psi^* (u_{-\omega_i}))(g)& =\Delta_{-\omega_i, -\omega_i}((n_c\oc)^{-1} m't) 
= t^{-\omega_i}\Delta_{w_0\omega_{i^*}, w_0\omega_{i^*}}((n_c\oc)^{-1} m')=
t^{-\omega_i}\Delta_{w_0\omega_{i^*}, w_0\omega_{i^*}}((n_c\oc)^{-1}) \\
&= t^{-\omega_i} \Delta_{\omega_i, \omega_i}
(n_c \oc) = t^{-\omega_i} \Delta_{\omega_i, \omega_i}(g) = f_i^{-1} u_{\omega_i}(g).
\end{align*}
Thus $\psi^*({\bf u}_{-\underline{\omega}}) = {\bf u}_{\underline{\omega}}\,{\bf f}^{-1} = 
{\bf u}_{\underline{\omega}}\,{\bf p}^{-U_{\bf 1}^{-1}}$. 
It also follows from $\psi(g)  = (\oc^{\, -1}n_c\oc)^{-1} t^c (\oc^{\, -1}) (t^{-1} m't)$ that 
\[
(\psi^*(f_i))(g) = f_i(\psi(g)) = t^{c\omega_i}= {\bf f}^{M_{\bf 1}e_i}.
\]
Thus $\psi^*({\bf f}) = {\bf f}^{M_{\bf 1}}$. It follows that 
$\psi^*({\bf p}) = \psi^*({\bf f}^{U_{\bf 1}})= {\bf f}^{M_{\bf 1} U_{\bf 1}} = {\bf f}^{-L_{\bf 1}}=
{\bf p}^{-U_{\bf 1}^{-1} L_{\bf 1}}$. 

Recall that $B_A = -L_{\bf 1} + U_{\bf 1}$. Taking $R = -U_{\bf 1}$ and $E = -U_{\bf 1}^{-1} {\bf L}_1$, one has
\[
RB_A + E = -U_{\bf 1}^{-1}(-L_{\bf 1} + U_{\bf 1}) - U_{\bf 1}^{-1}L_{\bf 1} = -I_n.
\]
By \ldref{:psi-extension}, $\psi: \Lcc \to \Lcc$ is a twisted 
automorphism of $\Lcc$ of DT-type, from which the statements on $\phi: \Acc \to \Acc$ and
$[\psi]:\Xcc \to \Xcc$ also follow. 
\end{proof}

\bre{:Weng}
{\rm
That $[\psi]$ is the DT transformation of $\Xcc$ also follows from 
\cite[Theorem 1.1]{Daping-DT-doubleBruhat} where DT transformations for general double Bruhat cells are studied.
\hfill $\diamond$
}
\ere

\bre{:psi-c-tau}
{\rm 
By \eqref{eq:Psi}, 
 $\psi^*: \CC[\Lcc] \to \CC[\Lcc]$ maps every
cluster variable to a cluster variable up to multiplication by a monomial of ${\bf p}$. 
Consider the bijection \cite{YZ:Lcc}
\[
\tau:\;\; \Pi(c) \longrightarrow \Pi(c), \;\;     \tau(c^m\omega_i) = 
    \begin{cases}
        c^{m+1}\omega_i & \text{ if } m \in [0, h(i;c)-1], \\
        \omega_{i^*} &\text{ if } m = h(i;c).
    \end{cases}
\] 
In addition to $\psi^*({\bf u}_{-\underline{\omega}}) = {\bf u}_{\underline{\omega}} \,{\bf p}^{-U_{\bf 1}^{-1}}$,
it is shown in \cite[Lemma 5.6.1 and Theorem 5.6.2]{Antoine:thesis} that 
\begin{equation}\label{eq:psi-gamma}
\psi^* u_{\gamma} =  u_{\tau(\gamma)}, \hs \forall\; 
\gamma \in \Pi(c)\backslash \{-\omega_1, \ldots, -\omega_n\}.
\end{equation}
In particular, $\psi^*(u_{\omega_i}) = u_{c\omega_i}$ for every $i \in [1,n]$.
Let again ${\bf h}$ be the Coxeter number of $\Phi$. 
 By \cite[Proposition 1.7]{YZ:Lcc}, one has
$h(i; c) + h(i^*; c) = {\bf h}$ for every $i \in [1, n]$. 
It follows that
$\tau^{{\bf h}+2} = {\rm Id}_{\Pi(c)}$, which gives another proof that $(\phi)^{{\bf h}+2} = {\rm Id}$.
\hfill $\diamond$
}
\ere


Turning now to  positive points, by \cite{FZ:double} the set $\Lccpos$ of positive points in $\Lcc$ defined by the Yang-Zelevinsky cluster structure coincides with 
$\Lcc \cap G_{\geq 0}$, where $G_{\geq 0}$ is the totally non-negative part of $G$ introduced by Lusztig 
\cite{Lusztig1994}, and $\Accpos = \Acc \cap G_{\geq 0}$. 
Set
\[
\psi_{ >0} = \psi|_{\Lccpos}:\;\; \Lccpos \longrightarrow \Lccpos
\hs \mbox{and} \hs {\phi}_{>0} = {\phi}|_{\Accpos}: \Accpos \to \Accpos,
\]
and let $[\psi]_{>0}=[\psi]|_{\Xccpos}: \Xccpos \to \Xccpos$. 
Recall again that 
$2 = d_2 \leq \cdots \leq d_n = {\bf h}$ are the degrees of the Weyl group $W$, and that with
$A = (a_{i, j} = \alpha_i^\vee(\alpha_j))$ as the Cartan matrix we have
\[
M_{\lambda} = -
\left(\begin{array}{ccccc} \lambda_1 & 0 & 0 & \cdots & 0\\
a_{2, 1} & \lambda_2 & 0 &\cdots & 0\\
a_{3,1} & a_{3,2} & \lambda_3 & \cdots & 0\\
\cdots & \cdots & \cdots & \cdots & \cdots\\
a_{n, 1} & a_{n, 2} & a_{n, 3} & \cdots & \lambda_n
\end{array}\right)\left(\begin{array}{ccccc} \lambda_1 & a_{1,2} & a_{1, 3} & \cdots & a_{1, n}\\
0 & \lambda_2 & a_{2,3} & \cdots & a_{2, n}\\
0 & 0 & \lambda_3 & \cdots & a_{3, n}\\
\cdots & \cdots & \cdots & \cdots & \cdots\\
0 & 0 & 0 & \cdots & \lambda_n\end{array}\right)^{-1}
\]
for $\lambda = (\lambda_1, \ldots, \lambda_n) \in (\CC^\times)^n$. Set again  ${\bf 1} = (1, \ldots, 1) \in (\CC^\times)^n$.
We now have the following reformulation of our main results in the paper in terms of the 
Fomin-Zelevinsky twist on $\Lcc$.

\bth{:Lie}
1) The map $\psi_{>0}$ has a unique fixed point
${\bf a} \in \Accpos \subset \Lccpos$, characterized by the 
property that $\Delta_{\omega_i, \omega_i}({\bf a}) = \Delta_{c\omega_i, c\omega_i}({\bf a})$ for all 
$i \in [1, n]$, or, equivalently, by 
$\Delta_{\gamma, \gamma}({\bf a}) = \Delta_{\tau \gamma, \tau \gamma}({\bf a})$ for every $\gamma \in \Pi(c)$;

2) Let $[{\bf a}]=p^{c, c^{-1}}({\bf a}) \in \Xcc$, and consider the three linear maps 
$d_{\bf a} (\psi_{>0}):  T_{\bf a} \left(\Lccpos\right)
\rightarrow  T_{\bf a} \left(\Lccpos\right)$,
\[
d_{\bf a} ({\phi}_{>0}):\;\; T_{\bf a} \left(\Accpos\right) \longrightarrow T_{\bf a}\left(\Accpos\right) \hs \mbox{and} \hs
d_{[{\bf a}]} ([\psi]_{>0}):\;\; T_{[{\bf a}]} \left(\Xccpos\right) \longrightarrow T_{[{\bf a}]}\left(\Xccpos\right).
\]
Let $\kappa = (\Delta_{\omega_1, \omega_1}({\bf a}), \ldots, \Delta_{\omega_n, \omega_n}({\bf a}))^A \in (\RR_{>0})^n$.
Then $d_{\bf a} ({\phi}_{>0})$ and $d_{[{\bf a}]} ([\psi]_{>0})$
have the same characteristic polynomial which is equal to 
\[
\calP_\Phi(x) = \det(x-M_{\kappa}) = \prod_{j = 1}^n \left(x - e^{\frac{2 \pi \sqrt{-1} }{{\bf h}+2}d_j}\right),
\]
and the characteristic polynomial of $d_{\bf a} (\psi_{>0})$ is equal to 
$\calP_\Phi(x) \calP_c(x)$, where 
\[
\calP_c(x) = \det(x-M_{\bf 1}) = \prod_{j = 1}^n \left(x - e^{\frac{2 \pi \sqrt{-1} }{{\bf h}}(d_j-1)}\right)
\] is the characteristic polynomial of the Coxeter element $c \in W$ as a linear operator on $\t^*$.
\eth

\begin{proof}
Suppose that ${\bf a} \in \Lccpos$ is fixed by $\psi_{>0}$, and let ${\bf q}= {\bf p}({\bf a})^{U_{\bf 1}}\in
(\RR_{>0})^n$. It then follows from \eqref{eq:FZtwist-u-p} that 
${\bf q}({\bf a}) ={\bf q}({\bf a})^{M_{\bf 1}}$, and so 
${\bf q}({\bf a})= {\bf 1}$ by \leref{:fixpt-frozen-1}. Thus ${\bf a} \in \calApos$, so ${\bf a}$ must be the only fixed point of $\phi$ in 
$\calApos$, whose characterization as given also follows from \eqref{eq:FZtwist-u-p}.

The statements on the characteristic polynomials of 
$d_{\bf a} ({\phi}_{ >0})$ and $d_{[{\bf a}]} ([\psi]_{>0})$ are proved in 
\prref{:matrix-realization-linearDT} and \thref{:expnts-conjecture}. 
The statement on the characteristic polynomial of $d_{\bf a} (\psi_{>0})$ follows from 
\thref{:same-char}, the fact that $\det(x+U_{\bf 1}^{-1}L_{\bf 1}) = \det(x+L_{\bf 1}U_{\bf 1}^{-1})
=\det(x-M_{\bf 1}) = \det(x-c)$ and that the exponents of $c$ with respect to 
${\bf h}$ are $d_1-1\leq d_2-1, \ldots, d_n-1$.
\end{proof}

\bex{:An-model}
Let $\Phi$ be of type $A_n$, with Cartan matrix $A = (a_{i,j})$ given by 
$a_{i, i}=2$ for all $i \in [1, n]$, $a_{i, j} = -1$ for all $|i-j|=1$ and $a_{i, j} =0$ otherwise. 
By \cite[Theorem 1.1]{YZ:Lcc}, $\Lcc = \{g({\bf q}, {\bf p}): {\bf q} \in \CC^{n+1}, {\bf p} \in (\CC^\times)^n, \det g({{\bf q}, {\bf p}}) = 1\}$, where
for ${\bf q} = (q_1, \ldots, q_{n+1}) \in \CC^{n+1}$ and ${\bf p} =(p_1, \ldots, p_n) \in (\CC^\times)^n$, 
\begin{equation}\label{eq:tridia-Acc}
   g= g({\bf q}, {\bf p}) = \begin{pmatrix}
        q_1 & p_1 & 0 & \cdots & 0 & 0\\
        1 & q_2 & p_2 & \cdots & 0 & 0\\
        0 & 1 & q_3 & \cdots & 0 & 0\\
        \vdots & \vdots & \vdots & \cdots & \vdots & \vdots \\
        0 & 0 & 0 & \cdots &q_n & p_n \\
        0 & 0 & 0 & \cdots & 1 & q_{n+1} \\
    \end{pmatrix}.
\end{equation}
Then the Fomin-Zelevinsky twist $\psi: \Lcc \to \Lcc$ is given by
\[
\psi(g((q_1, \ldots, q_{n+1}, p_1, \ldots, p_n)) 
= g(q_2, \ldots, q_{n+1}, q', \,p_2, \,\ldots,\, p_n, \,p_1^{-1}p_2^{-1}\cdots p_n^{-1}),
\]
where $q' \in \CC$ is such that 
$\det g(q_2, \ldots, q_{n+1}, q', p_2, \ldots, p_n, p_1^{-1}p_2^{-1}\cdots p_n^{-1}) = 1$.

Let $g(q_1,\ldots, q_{n+1}) = g(q_1, \ldots, q_{n+1}, 1, \ldots, 1)$.
Then $\Acc=\{g = g(q_1, \dots, q_{n+1}): \det(g) = 1\}$. 
For $j \in [1, n+1]$, let $\varphi_j := \varphi_j(q_1, \ldots, q_{n+1})$ be the $j^{\rm th}$ principal minor of 
$g(q_1, \ldots, q_{n+1})$.
 Then 
\begin{equation}\label{eq:phi-j}
\varphi_j = q_j \varphi_{j-1} - \varphi_{j-2}, \hs j \in [1, n+1],
\end{equation}
where $\varphi_0 = 1$ and $\varphi_{-1} = 0$. Identifying $\Acc$ with the sub-variety $Q$ of $\CC^{n+1}$ given by the equation
\[
q_{n+1}\varphi_n -\varphi_{n-1}=1,
\] 
the Fomin-Zelevinsky twist on $\Acc$ now becomes 
\begin{equation}\label{eq:psi-An}
\phi: \;\; Q \longrightarrow Q, \;\; \phi(q_1, \ldots,q_n,  q_{n+1}) = (q_2, \ldots, q_{n+1}, 
\varphi_n).
\end{equation}
One checks directly that $\phi$ is periodic with period $n+3$. 
Fixed points of $\phi$ correspond to all $q \in \CC$ such that 
$\det g(q, \ldots, q) = 1$, equivalently,  $v_{n+1}(q) = 1$, where $(v_j)_{j \in \mathbb{N}}$ is the sequence of {\it Vieta polynomials} defined by 
\[
v_0 = 1, v_1 = q \quad \text{and} \quad v_j = q v_{j-1} - v_{j-2} \quad \text{for } j \geq 2,
\]
 and they are related to the Chebyshev polynomials 
$(U_j)_{j \in \mathbb{N}}$ of the second kind  by $v_j(q) = U_j(\frac{q}{2})$; see \cite{Koshy:Vieta}. One checks directly that solutions to $v_{n+1}(q)=1$ are given by $q = 2\cos \theta$, where 
\[
    \theta = \frac{2k\pi}{n+1}, \, k \in \left\{1,2,\ldots, \left\lfloor \frac{n+1}{2} \right\rfloor \right\} \quad \text{or} \quad  \theta = \frac{(2k-1)\pi}{n+3}, \, k \in \left\{ 1,2,\ldots , \left\lfloor \frac{n+2}{2} \right\rfloor\right\} .
\]
By \thref{:Lie}, there is exactly {\it one} solution  $q \in \CC$ such that  
$g(q, \ldots, q) \in SL(n+1, \CC)$ is a totally non-negative matrix,  i.e. with all minors non-negative.
One checks directly that 
$q = 2\cos \left(\frac{\pi}{n+3}\right)$.
In $\S$\ref{ss:appendix:An-Cn} we will use the
linearization of $\phi$ at $g(q, \ldots, q)$ to prove \thref{:expnts-conjecture} for type $A_n$.
$\hfill \diamond $
\eex

\appendix

\section{Fixed points in each Cartan type and the polynomials $\calD_{X_n}(x)$}\label{s:appendixA}

\setcounter{equation}{0}\renewcommand\theequation{A\arabic{equation}}

\begin{table}[h]
    \centering
\[
\begin{array}{|c|c|c|}
    \hline
    \text{Cartan type} & \text{Dynkin diagram}&  \text{Positive solution $(a_1, \ldots, a_n)$ of \eqref{eq:fixpt-eq-A}} \\\hline
     A_n& \includegraphics[scale=0.6]{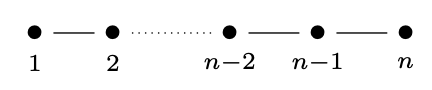} &\left(\frac{\sin(2\theta_n)}{\sin(\theta_n)}, \;
\frac{\sin(3\theta_n)}{\sin(\theta_n)}, \;\ldots, \; 
\frac{\sin(n+1)\theta_n}{\sin(\theta_n)}\right)\\[3ex]
B_n & \includegraphics[scale=0.6]{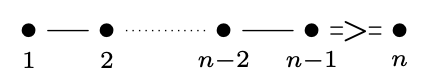} &\left(2, \; 3, \ldots, n, \sqrt{n+1} \right) \\[3ex]
     C_n & \includegraphics[scale=0.6]{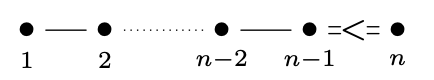}& \left(\frac{\sin(2\theta_{2n-1})}{\sin(\theta_{2n-1})}, \;
\frac{\sin(3\theta_{2n-1})}{\sin(\theta_{2n-1})}, \;\ldots, \; 
\frac{\sin(n+1)\theta_{2n-1}}{\sin(\theta_{2n-1})}\right) \\[3ex]
     D_n & \includegraphics[scale=0.6]{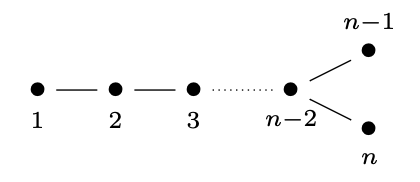} &\left(2, \; 3, \ldots, n-1, \sqrt{n}, \sqrt{n} \right) \\[3ex]
     E_6 &\includegraphics[scale = 0.6]{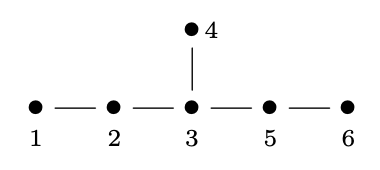} &(\varepsilon, \;\; \varepsilon^2-1, \;\; 2\varepsilon^2-\varepsilon-1, \;\; 
\varepsilon(\varepsilon-1), \;\; \varepsilon^2-1, \;\;\varepsilon) \\[3ex]
     E_7 &\includegraphics[scale =0.6]{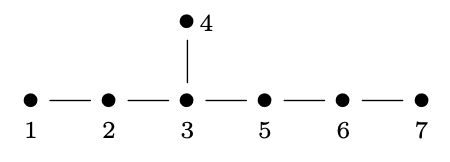}& \Delta_1\\[3ex]
     E_8 &\includegraphics[scale=0.6]{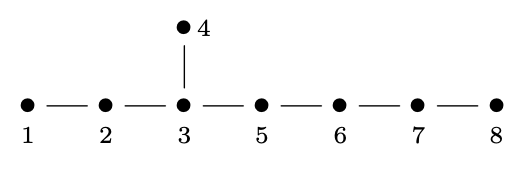}& \Delta_2\\[3ex]
     F_4 &\includegraphics[scale=0.6]{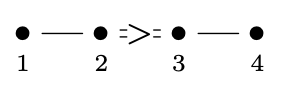}& (\varepsilon, \;\;\varepsilon^2-1, \;\; 2\varepsilon^2-\varepsilon-1,\;\;
\varepsilon(\varepsilon-1)) \\[3ex]
     G_2 & \includegraphics[scale=0.6]{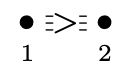} &(3,2) \\
     \hline
\end{array}
\]
$\theta_n = \frac{\pi}{n+3}$,  $\varepsilon = 2.2470...$ is the
unique real root in $(2, \infty)$ of the equation
$\varepsilon^3 - 2 \varepsilon^2 - \varepsilon + 1=0$, and
\begin{align*}
    \Delta_1 &= ((3+\sqrt{5})/2, \; (5+3\sqrt{5})/2,\; 6+3\sqrt{5},\; \sqrt{7 +3\sqrt{5}},\; 2\sqrt{7 +3\sqrt{5}},\;  2+\sqrt{5}, \;\sqrt{3 +\sqrt{5}}), \\
    \Delta_2 &= (2+\sqrt{2}, \; 5+4 \sqrt{2},\; 16+12 \sqrt{2}, \;3+2\sqrt{2}, \; 9+6\sqrt{2}, \; 5+3\sqrt{2}, \;
 2 + 2\sqrt{2}, \; 1+\sqrt{2}).
\end{align*}
    \caption{Coordinates of fixed points of the DT transformations}
    \label{table:A}
\end{table}

\begin{table}[h]
    \centering
\[
\begin{array}{|c|c|c|}
    \hline
    \text{Cartan type} &
d_1, \ldots, d_n &  {\bf h}  \\
\hline
A_n &2, 3, \ldots, n+1 & n+1 \\
B_n, C_n & 2, 4, 6, \ldots, 2n & 2n\\
D_n &2, 4,6, \ldots, 2n-2, n & n\\
E_6 &  2,5,6,8,9,12 & 12\\ 
E_7 & 2,6,8,10,12,14,18 & 18\\
E_8 & 2,8,12,14,18,20,24,30 & 30\\
F_4 & 2,6,8,12 & 12 \\
G_2 & 2, 6 & 6\\
\hline
\end{array}
\]
\caption{Degrees and Coxeter number of the Weyl group}
    \label{table:d}
\end{table}

For each indecomposable finite root system $\Phi$, 
we give in Table \ref{table:A} the Dynkin diagram of $\Phi$ with a labeling of 
its vertices, which determines the Cartan matrix $A = (a_{i, j})$, and the corresponding unique solution
${\bf a} = (a_1, \ldots, a_n)$  to the equations in \eqref{eq:fixpt-eq-A}. By \prref{:ab-equations}, 
${\bf b} = (b_1 , \ldots , b_n) = (a_1^2-1, \, \ldots, \, a_n^2-1)$
is the unique positive solution to the equations in \eqref{eq:fixpt-eq-X}.

For the convenience of the reader, we recall in  Table \ref{table:d} the 
explicit values of the Coxeter number ${\bf h}$ of $\Phi$ and of the degrees $2 =d_1  \leq d_2 \leq \cdots \leq d_n ={\bf h}$ of the Weyl group of $\Phi$
from \cite[Table 1, p.59]{HumphreysCoxGrp}).
Recall that 
\begin{equation*}
    \calD_\Phi= \prod_{j = 1}^n \left(x - e^{\frac{2\pi \sqrt{-1} \,d_j}{{\bf h}+2}}\right).
\end{equation*}
If $\Phi$ is of type $X_n$, we write $\calD_{X_n}(x) = \calD_\Phi(x)$.

\ble{:D-poly} The polynomials
$\calD_{X_n}(x)$ are given as follows:
\begin{align}\label{eq:D-An}
\calD_{A_n}(x) &= 
\sum_{j=0}^{n} \frac{\sin ((j+1)\theta_n) \sin ((j+2)\theta_n)}{\sin (\theta_n) \sin (2\theta_n)}x^j;\hs n \geq 1, \;\;
\theta_n =\frac{\pi}{n+3};\\
\label{eq:D-Bn}
\calD_{B_n}(x) &= \calD_{C_n}(x) = 1 + x + \ldots + x^n, \quad n \geq 2, \\
\label{eq:D-Dn}
\calD_{D_n}(x) &= 1 + 2x + 2x^2 + \ldots + 2x^{n-1} + x^n, \quad n \geq 4,\\
\label{eq:D-G2}
\calD_{G_2}(x) &= 1 + x^2,\\
\label{eq:D-F4}
\calD_{F_4}(x) &= x^4 + \alpha_1x^3 + \alpha_2x^2 + \alpha_1x + 1, \\
\label{eq:D-E6}
\calD_{E_6}(x) &= x^6 + \beta_1x^5 + \beta_2x^4 + \beta_3x^3 + \beta_2x^2+\beta_1x+ 1, \\  
\label{eq:D-E7}
\calD_{E_7}(x) &= x^7 + \tau x^6 + x^5+x^2+\tau x+ 1, \quad \tau = \frac{1 + \sqrt{5}}{2},\\
\label{eq:D-E8}    
\calD_{E_8}(x) &= x^8 + \gamma_1x^7 + \gamma_2x_6 +\gamma_3x^5+ \gamma_4x^4 + \gamma_3x_3 +\gamma_2x^2+\gamma_1x+ 1,
\end{align}
where, with $\omega = e^{\frac{\pi \sqrt{-1}}{7}}$,
\begin{align*}
    \alpha_1 &= -(\omega^2+\omega^{-2}+\omega^6+\omega^{-6}),\hs \alpha_2 = 2+\omega^4+\omega^{-4}+\omega^6+\omega^{-6},\\
    \beta_1 &= \omega +\omega^{-1}, \hs \beta_2 = 1 + \omega^3 +\omega^{-3}, \hs \beta_3 = \omega^2+\omega^{-2}+\omega^4+\omega^{-4},\\
    \gamma_1 &= \sqrt{2}, \hs \gamma_2 = 2 - \sqrt{2}, \hs \gamma_3 = \sqrt{2}-2, \hs \gamma_4 = 2-2\sqrt{2}.
\end{align*}
\ele

\begin{proof}
We only prove the case of $A_n$, others cases being the results of direct calculations.
Let $\zeta_n = e^{\frac{2\pi \sqrt{-1}}{n+3}}$, so that $\calD_{A_n}(x) = \prod_{j=2}^{n+1}(x-\zeta_n^j)$. As
$(x-\zeta_n) (x-\zeta_n^{-1}) = x^2 - 2\cos(2\theta_n)x + 1$, one has
\[
 (x^2 - 2\cos(2\theta_n)x + 1) \calD_{A_n}(x) = \frac{x^{n+3}-1}{x-1} = 1 + x + x^2 + \ldots + x^{n+2}.
\]
Write $\calD_{A_n}(x) = \sum_{j=0}^n c^{(n)}_j x^j$ and set $c^{(n)}_{-1} = c^{(n)}_{-1} =0$.  
Then $\{c^{(n)}_j\}_{j \in [0, n]}$ are recursively determined by
\begin{equation}\label{eq:cnj}
c^{(n)}_j -2 \cos(2\theta_n) c^{(n)}_{j-1} + c^{(n)}_{j-2} = 1, \hs j \in [0, n].
\end{equation}
One checks directly that the only solution to \eqref{eq:cnj} is
$c^{(n)}_j = \frac{\sin ((j+1)\theta_n) \sin ((j+2)\theta_n)}{\sin (\theta_n) \sin (2\theta_n)}$ for $j \in [0, n]$.
\end{proof}

\section{Proof of \thref{:expnts-conjecture} (Theorem B)}\label{s:appendixB}
\setcounter{equation}{0}\renewcommand\theequation{B\arabic{equation}}

For a root system $\Phi$ of Cartan type $X_n$, write $\calP_{X_n}(x)$ for the 
characteristic polynomial of the linear maps $\dtA$ and $\dtX$ 
for the cluster ensemble $(\calA, \calX, p)$ of type $\Phi$. In this Appendix, we prove that 
\begin{equation}\label{eq:key}
\calP_{X_n}(x) = \calD_{X_n}(x),
\end{equation}
where $\calD_{X_n}(x)$ is given in \leref{:D-poly}.
By \prref{:matrix-realization-linearDT}, we have 
\begin{equation}\label{eq:key-det}
    \calP_{X_n}(x) = \det(xI - M_{\bf \kappa}) = \frac{1}{\prod_{i=1}^n \ka_i}\det(U_{\bf \ka}x + L_{\bf \ka}),
\end{equation}
where, with ${\bf a} = (a_1,\ldots, a_n)$ given in Table \ref{table:A} and $A$ denoting the corresponding Cartan matrix,
\[
\kappa = (\kappa_1, \ldots, \kappa_n) =\left(\frac{a_1^2}{a_1^2-1}, \ldots, \frac{a_n^2}{a_n^2-1}\right) = 
{\bf a}^A,
\]
and $U_\ka$ and $L_\ka$ are given in \eqref{eq:LU-intro}.  When $X_n$ can be obtained by folding
from $Y_m$ of simply laced type, instead of applying  folding in the 
context of cluster algebras \cite{FSP:unfolding, Huang-Li:unfolding}, we use \eqref{eq:key-det} directly to write $\calP_{Y_m}(x)$ as a product with $\calP_{X_n}(x)$ as a factor and simplify some of the case-by-case proofs. 

\subsection{Type $A_n$ and Type $C_n$}\label{ss:appendix:An-Cn} We first consider type $A_n$ for $n \geq 1$. For $j \in [1, n]$, let $\varphi_j(q_1, \ldots, q_j)$ be the
Euler continuant given in \eqref{eq:phi-j}, and let $q = 2\cos\theta_n$, where $\theta_n =\frac{\pi}{n+3}$. 
By  \exref{:An-model}, 
the co-dimensional one sub-variety
$Q$ of $\CC^{n+1}$  defined by
\begin{equation}\label{eq:q-phi}
q_{n+1} \varphi_{n}(q_1, \ldots, q_n) - \varphi_{n-1}(q_1, \ldots, q_{n-2}) = 1,
\hs (q_1, \ldots, q_{n+1}) \in \CC^{n+1},
\end{equation}
is an affine model of $\calA$, using which  $\DTA: \calA \to \calA$ 
becomes the bi-regular isomorphism 
\[
\phi:\;\; Q \longrightarrow Q, \;\; \phi(q_1, \,q_2,\, \ldots,\, 
q_{n+1}) = (q_2,\, \ldots,\, q_n, \, \varphi_n(q_1, \ldots, q_n)),
\]
and ${\bf q} = (q, q, \ldots, q) \in Q$ is the unique positive fixed point of $\phi$. Thus 
$\calP_{A_n}(x)$ is the characteristic polynomial of the linear map
$(d_{\bf q}\phi)^*: T_{\bf q}^*Q \to T_{\bf q}^*Q$. Regard $q_1, \ldots, q_n$ as regular functions on $Q$ and
consider the basis 
$(\delta_1, \ldots, \delta_n) = (d_{\bf q} q_1, \ldots, d_{\bf q} q_n)$ of $T^*_{\bf q} Q$. We then have
\begin{equation}\label{eq:d-delta}
(d_{\bf q} \phi)^* (\delta_1, \ldots, \delta_n) = (\delta_2, \,\ldots, \,\delta_{n},\,d_{\bf q} q_{n+1}).
\end{equation}
Thus $\calP_{A_n}(x) = c_1 + c_2x + \cdots + c_nx^{n-1} + x^n$ if $d_{\bf q} q_{n+1}=-(c_1\delta_1+ \cdots + c_n \delta_n)$.
Set ${\bf a}_{A_n} = (a_1, \ldots, a_{2n-1})$.
One 
computes directly that 
\begin{equation}\label{eq:a-n}
{\bf a}_{A_n}= (\varphi_1(q), \; \varphi_2(q, q), \; \ldots, \;\varphi_n(q, \ldots, q)) = \left(\frac{\sin(2\theta_n)}{\sin(\theta_n)}, \;
\frac{\sin(3\theta_n)}{\sin(\theta_n)}, \;\ldots, \; 
\frac{\sin(n+1)\theta_n}{\sin(\theta_n)}\right).
\end{equation}
For $j \in [1, n]$, regarding $\varphi_j =\varphi_j(q_1, \ldots, q_j)$ as regular functions on $Q$, 
setting $\sigma_j = d_{\bf q}\varphi_j \in T_{\bf q}^* Q$, and differentiating \eqref{eq:q-phi} at ${\bf q}$, we then have $-qd_{\bf q} q_{n+1} = q\sigma_n -\sigma_{n-1}$. Using 
$\varphi_j = q_j \varphi_{j-1}-\varphi_{j-2}$ for $j \in [1, n]$. 
where $\varphi_0=1$ and $\varphi_{-1} = 0$ by definition, and setting $\sigma_0 = \sigma_{-1} = 0$, we get
$\sigma_j = a_{j-1} \delta_j + q\sigma_{j-1} -\sigma_{j-2}$ for $j \in [1, n]$.,
An inductive argument using  the fact that $a_j = a_{n+1-j}$ for $j \in [1, n]$ then gives 
\[
-qd_{\bf q} q_{n+1} = a_1a_{n-1} \delta_n + (q^2-1)\sigma_{n-1}-q\sigma_{n-2} =
a_1a_{n-1} \delta_n + a_2a_{n-1} \delta_{n-1} + \cdots + a_{n-1}a_1 \delta_1.
\]
Comparing with \eqref{eq:D-An}, we see that $\calP_{A_n}(x) = \calD_{A_n}(x)$.

We now turn to type $C_n$ for $n \geq 3$. Note first that $C_n$ can be obtained from $A_{2n-1}$ by folding. In particular, $C_n$ and $A_{2n-1}$ have the same Coxeter number ${\bf h} = 2n$. Set $\zeta = e^{\frac{2 \pi \sqrt{-1}}{2n+2}} =
e^{\frac{\pi \sqrt{-1}}{n+1}}$. Then  
\[
\calP_{A_{2n-1}}(x) = \calD_{A_{2n-1}}(x) = \prod_{j=2}^{2n} (x-\zeta^j) 
=\prod_{j=1}^n (x-\zeta^{2j})\prod_{j=2}^{n}(x-\zeta^{2j-1}) = \calD_{C_n}(x) \prod_{j=2}^{n}(x-\zeta^{2j-1}).
\]
By \leref{:DT-periodic}, $\DTA$ of type $C_n$ has period $n+1$, so all the (complex)
roots of $\calP_{C_n}(x)$ are even powers of $\zeta$. Since $\calP_{C_n}(x)$ and $\calD_{C_n}(x)$ have the same degree, to prove that $\calP_{C_n}(x) = \calD_{C_n}(x)$, 
it suffices to prove that $\calP_{C_n}(x)$ is a factor of $\calP_{A_{2n-1}}(x)$ in $\CC[x]$.
To this end, note the linear algebra fact that 
\begin{equation}\label{eq:det-det}
\det\begin{pmatrix}
    \alpha_1 & \beta_1 & 0 &\cdots  &0 &0\\
    \gamma_1       & \alpha_2 & \beta_2 &\cdots &0&0\\
    0      & \gamma_2 & \alpha_3  &\cdots&0&0\\
    \cdots & \cdots & \cdots & \cdots  & \cdots& \cdots \\
     0 & 0 &0  &\cdots   & \alpha_{m-1}& \beta_{m-1} \\
    0 & 0 & 0 & \cdots & \gamma_{m-1} &\alpha_m
 \end{pmatrix} 
  = \det \begin{pmatrix}
    \alpha_m & \beta_{m-1} & 0 &\cdots  &0 &0\\
    \gamma_{m-1}       & \alpha_{m-1} & \beta_{m-2} &\cdots &0&0\\
    0      & \gamma_{m-2} & \alpha_{m-2}  &\cdots&0&0\\
    \cdots & \cdots & \cdots & \cdots  & \cdots& \cdots \\
     0 & 0 &0  &\cdots   & \alpha_{2}& \beta_{1} \\
    0 & 0 & 0 & \cdots & \gamma_{1} &\alpha_1
 \end{pmatrix} 
\end{equation}
for $m \geq 2$ and any elements
$\alpha_1, \ldots, \alpha_m, \beta_1, \ldots, \beta_{m-1}, \gamma_1, \ldots, \gamma_{m-1}$ in any commutative ring
(the second matrix is the transpose of a conjugate of the first by a permutation matrix). 
Let now
${\bf a}_{A_{2n-1}} = (a_1, \ldots, a_{2n-1})$ be as in \eqref{eq:a-n} and let 
$\kappa_j = \frac{a_j^2}{a_j^2-1}$ for $j \in [1, 2n-1]$. By \eqref{eq:key-det}, 
$\calP_{A_{2n-1}} = \frac{1}{\ka_1\ka_2 \cdots \ka_{2n-1}}\det M$, where 
\[
M=\begin{pmatrix}
    \ka_1(x+1) & -x & 0 & 0 &\cdots &  0 & 0 &0\\
    -1       & \ka_2(x+1) & -x &  0 &\cdots &0 & 0&0\\
    \cdots & \cdots & \cdots & \cdots & \cdots & \cdots& \cdots&\cdots \\
     0 & 0 &0  & 0 &\cdots  &  -1 & \ka_{2n-2}(x+1)& -x \\
    0 & 0 & 0 & 0 & \cdots & 0 & -1 &\ka_{2n-1}(x+1)
 \end{pmatrix}.
\]
Let $M_1$ be the principal $(n-1) \times (n-1)$ sub-matrix of $M$ and 
$M_2$ the principal $(n-2) \times (n-2)$ sub-matrix of $M$. 
Expanding $\det(M)$ along the $n^{\rm th}$ row of $M$ 
and using \eqref{eq:det-det},  it is straightforward to show that 
\begin{align*}
\calP_{A_{2n-1}} &= \frac{1}{\ka_1\ka_2 \cdots \ka_{2n-1}}\left(\ka_n(x+1)\det (M_1)^2-2x\det(M_1)\det(M_2)\right)
=\frac{1}{\ka_{n+1} \cdots \ka_{2n-1}}\det(M_1)\calP_{C_n}(x),
\end{align*}
where we used the fact that $\calP_{C_n} =\frac{1}{\ka_1\cdots \ka_n}( \ka_n (x+1)\det(M_1) - 2x \det(M_2))$ 
by \eqref{eq:key-det} and 
Table \ref{table:A}.
We have thus proved that $\calP_{C_n}(x)$ is a factor of $\calP_{A_{2n-1}}$. Thus
$\calP_{C_n}(x) = \calD_{C_n}(x)$. 
We also proved that 
\[
 \frac{1}{\ka_1 \cdots \ka_{n-1}}\det(M_1) = \prod_{j=2}^{n}(x-\zeta^{2j-1}).
\]

\subsection{Type $B_n$ and Type $D_n$} A direct computation gives 
$\calP_{B_2}(x) = 1 + x + x^2 = \calD_{B_2}(x)$. Let $n \geq 4$. By \eqref{eq:key-det} and Table \ref{table:A} for type $D_n$, we have 
\begin{equation}\label{eq:Dn-char-poly}
    \calP_{D_n}(x)  =\frac{n-1}{2n} \det \left(
\begin{array}{ccccccc} \frac{4(x+1)}{3} & -x  & \cdots & 0 & 0 & 0 & 0\\
-1 & \frac{9(x+1)}{8} & \cdots & 0 & 0 & 0 & 0\\
\vdots & \vdots  & \vdots & \vdots & \vdots & \vdots & \vdots \vspace{.03in} \\
0 & 0 & \cdots &   \frac{(n-2)^2(x+1)}{(n-1)(n-3)}  & -x & 0 & 0 \\
0 & 0 & \cdots & -1 &\frac{(n-1)^2(x+1)}{n(n-2)} & -x & -x \\
0 & 0  & \cdots & 0 & -1 & \frac{n(x+1)}{n-1} &  0  \\
0 & 0  & \cdots & 0 & -1 & 0 & \frac{n(x+1)}{n-1}\end{array}\right).
\end{equation}
Subtracting $n^{\rm th}$ column from $(n-1)^{{\rm th}}$ column, one has 
\begin{align}\nonumber
\calP_{D_n}(x)&=
\frac{x+1}{2} \det \left(
\begin{array}{ccccccc} \frac{4(x+1)}{3} & -x  & \cdots & 0 & 0 & 0 & 0\\
-1 & \frac{9(x+1)}{8} & \cdots & 0 & 0 & 0 & 0\\
\vdots & \vdots  & \vdots & \vdots & \vdots & \vdots & \vdots  \\
0 & 0 & \cdots &\frac{(n-2)^2(x+1)}{(n-1)(n-3)}   & -x & 0 & 0 \\
0 & 0 & \cdots & -1 &\frac{(n-1)^2(x+1)}{n(n-2)} & 0 & -x \\
0 & 0 & \cdots & 0 & -1 & 1 &  0  \\
0 & 0 & \cdots & 0 & -1 & -1 & \frac{n(x+1)}{n-1}\end{array}\right)\\
\label{eq:P-Dn}& = 
\frac{x+1}{2} \det\left(\begin{array}{cccccc} \frac{4(x+1)}{3} & -x & \cdots & 0 & 0 & 0\\
-1 & \frac{9(x+1)}{8} & \cdots & 0 & 0  & 0\\
\vdots & \vdots & \vdots & \vdots & \vdots & \vdots \vspace{.03in} \\
0 & 0 & \cdots & \frac{(n-2)^2(x+1)}{(n-1)(n-3)}  & -x  & 0 \\
0 & 0 & \cdots & -1 &\frac{(n-1)^2(x+1)}{n(n-2)} & -x \\
0 & 0 & \cdots & 0 & -2  & \frac{n(x+1)}{n-1}\end{array}\right).
\end{align}
Note that the determinant on the right hand side of \eqref{eq:P-Dn} is that of an $(n-1) \times (n-1)$ matrix. 
By \eqref{eq:key-det} and Table \ref{table:A} for type $B_{n-1}$, we get 
\begin{equation}\label{eq:P-Dn-Bn}
\calP_{D_n}(x)  = (x+1)\calP_{B_{n-1}}(x).
\end{equation}
Comparing with $\calD_{D_n}(x)$ and $\calD_{B_{n-1}}(x)$ in \leref{:D-poly}, it suffices to prove that 
$\calP_{B_{n-1}}(x) = 1 + x  + \cdots + x^{n-1}$. 

Expanding the determinant on the right hand side of in \eqref{eq:P-Dn} along the last row, we have
\begin{equation}\label{eq:P-B-n}
\calP_{B_{n-1}}(x) = \frac{n(x+1)}{2(n-1)} \Delta_{n-2}(x) - x \Delta_{n-3}(x),
\end{equation}
where for $m\geq 1$, $\Delta_m(x)$ is the determinant of an $m \times m$ matrix
\[
\Delta_m(x) = \det \begin{pmatrix}\frac{4(x+1)}{3} & -x & \cdots & 0 & 0\\
-1 & \frac{9(x+1)}{8} & \cdots & 0& 0\\
\vdots & \vdots & \vdots & \vdots& \vdots\\
0 & 0 & \cdots & \frac{m^2(x+1)}{(m-1)(m+1)}& -x\\
0 & 0 & \cdots & -1 & \frac{(m+1)^2(x+1)}{m(m+2)} 
\end{pmatrix}.
\]
We also set $\Delta_1(x) = 1$. Applying \eqref{eq:P-B-n} to $B_{n-2}$ and using 
\begin{equation}\label{eq:Delta-m}
\Delta_m = \frac{(m+1)^2(x+1)}{m(m+2)}\Delta_{m-1}(x)-x\Delta_{m-2}(x), \hs m \geq 2,
\end{equation}
we get $\calP_{B_{n-1}}(x)-\calP_{B_{n-2}}(x) = xR_{n-2}$, where
\[
R_{n-2}  =\frac{(n-1)x-(n-3)}{2(n-2)}\Delta_{n-3}(x) - \frac{nx-(n-2)}{2(n-1)} \Delta_{n-4}(x).
\]
Using again \eqref{eq:Delta-m} for $m = n-3$, one sees that $R_{n-2} = xR_{n-3} = \cdots = x^{n-4}R_2 = x^{n-2}$.
Thus
\[
\calP_{B_{n-1}}(x) = x^{n-1} + \calP_{B_{n-2}}(x)  = \cdots = x^{n-1} + \cdots + x + 1.
\]
Thus finishes the proof that $\calP_{B_n}(x) = \calD_{B_n}(x)$ for $n \geq 2$ and 
$\calP_{D_n}(x) = \calD_{D_n}(x)$ for $n \geq 4$.

\subsection{Type $G_2$}
We have
\[
 \calP_{G_2}(x) = \frac{2}{3}\det \begin{pmatrix}
     \frac{9}{8}(x+1) & -x \\
     -3 & \frac{4}{3}(x+1)
 \end{pmatrix} = x^2+1 = \calD_{G_2}(x).
\]

\subsection{Type $E_6$ and Type $F_4$} Let $(a_1, \ldots, a_6)$ be as in Table \ref{table:A} for $E_6$ and
$\ka_i = \frac{a_i^2}{a_i^2-1}$ for $i \in [1, 6]$. Then 
\begin{align*}
\calP_{E_6}(x) &= \frac{1}{\ka_1\ka_2\cdots \ka_6}\det \begin{pmatrix}
    \ka_1(x+1) & -x & 0&0&0&0 \\
    -1& \ka_2(x+1) & -x &0&0&0\\
    0&-1&\ka_3(x+1)&-x&-x&0\\
    0&0&-1&\ka_4(x+1)&0&0\\
    0&0&-1&0&\ka_5(x+1)&-x\\
    0&0&0&0&-1&\ka_6(x+1)
 \end{pmatrix}.
\end{align*}
Expanding the determinant along the $4^{\rm th}$ column and setting $f_2 = \frac{1}{\ka_1\ka_2}(\ka_1\ka_2(x+1)^2-x)$, we get
\begin{align*}
\calP_{E_6}(x) &= \frac{x}{\ka_3\ka_4}f_2^2 +\frac{x+1}{\ka_1\ka_2\ka_3\ka_5\ka_6}\det
\begin{pmatrix}
    \ka_1(x+1) & -x & 0&0&0 \\
    -1& \ka_2(x+1) & -x &0&0\\
    0&-1&\ka_3(x+1)&-x&0\\
    0&0&-1&\ka_5(x+1)&-x\\
    0&0&0&-1&\ka_6(x+1)
 \end{pmatrix}\\
 &  = \left(-\frac{x}{\ka_3\ka_4}f_2 +(x+1)^2f_2
 -\frac{2x(x+1)^2}{\ka_2\ka_3}\right)f_2.
 \end{align*}
On the other hand, by \eqref{eq:key-det} and Table \ref{table:A} for type $F_4$, we have 
\begin{align*}
\calP_{F_4}(x) &= \frac{1}{\ka_1\ka_2\ka_3\ka_4}\det \begin{pmatrix}
     \ka_1 (x+1) & -x & 0&0 \\
     -1 & \ka_2(x+1) & -2x & 0\\
     0 & -1 & \ka_3(x+1) & -x \\
     0&0&-1&\ka_4(x+1)
 \end{pmatrix} \\& = -\frac{x}{\ka_3\ka_4}f_2 +(x+1)^2f_2
 -\frac{2x(x+1)^2}{\ka_2\ka_3}.
 \end{align*}
One thus has $\calP_{E_6}(x) =\calP_{F_4}(x) f_2$. Let $\omega = e^{\frac{\pi \sqrt{-1}}{7}}$.
With $\varepsilon$ given in Table \ref{table:A}, one has 
\[
\varepsilon = 2\cos\frac{\pi}{7} + 2\cos\frac{3\pi}{7} = \omega + \omega^{-1}+\omega^3 + \omega^{-3}.
\]
For $\alpha_1, \alpha_2$ given in \leref{:D-poly}, one has 
$\alpha_1 = -2\varepsilon + \varepsilon^2$ and  $\alpha_2 = 2-\varepsilon$. 
Using the explicit values of $\ka_1, \ka_2, \ka_3, \ka_4$ in terms of $\varepsilon$, one then checks that 
$\calP_{F_4}(x) = \calD_{F_4}(x)$ and $\calP_{E_6}(x) = \calD_{E_6}(x)$.

\subsection{Type $E_7$} With $(a_1,\ldots, a_7)$ given as in Table \ref{table:A} for Type $E_7$
and $\ka_j = \frac{a_j^2}{a_j^2-1}$ for $j \in [1, 7]$, we use \eqref{eq:key-det} to compute
$\calP_{E_7}(x)$ to get
\[
\calP_{E_7}(x)= x^7 + \lambda_1x^6 + \lambda_2x^5 + \lambda_3x^4 + \lambda_3x^3+\lambda_2x^2+\lambda_1x+ 1
\]
where $\lambda_1 = 7 -\frac{c_1}{\ka}$, $\lambda_2 = 21 -5 \frac{c_1}{\ka} +\frac{c_2}{\ka}$,
and  $\lambda_3 = 35-10\frac{c_1}{\ka}+3\frac{c_2}{\ka}-\frac{c_3}{\ka}$ with
$\ka = \ka_1\ka_2\ka_3\ka_4\ka_5\ka_6\ka_7$, and 
\begin{align*}
    c_1 &=\ka_1 \ka_2 \ka_3 \ka_4 \ka_5 + \ka_1 \ka_2 \ka_3 \ka_4 \ka_7 + \ka_1 \ka_2 \ka_4 \ka_6 \ka_7 + \ka_1 \ka_2 \ka_5 \ka_6 \ka_7 + 
 \ka_1 \ka_4 \ka_5 \ka_6 \ka_7 + \ka_3 \ka_4  \ka_5  \ka_6  \ka_7,\\
    c_2 &= \ka_1 \ka_2 \ka_4 + \ka_1 \ka_2 \ka_5 + \ka_1 \ka_4 \ka_5 + \ka_3 \ka_4 \ka_5  + 
 \ka_1 \ka_2 \ka_7 + \ka_1 \ka_4 \ka_7 + \ka_3 \ka_4 \ka_7  + \ka_4 \ka_6 \ka_7  + \ka_5 \ka_6 \ka_7, \\
 c_3 &= \ka_4  +\ka_5  +\ka_7.
\end{align*}
Using the explicit values of $(a_1, \ldots, a_7)$ in  Table \ref{table:A}, we get 
 $\lambda_1=\tau, \lambda_2 = 1$ and $\lambda_3 = 0$. Thus $\calP_{E_7}(x) = \calD_{E_7}(x)$.

\subsection{Type $E_8$}
With $(a_1,\ldots, a_8)$ given as in Table \ref{table:A} for Type $E_8$
and $\ka_j = \frac{a_j^2}{a_j^2-1}$ for $j \in [1, 8]$, we use \eqref{eq:key-det} to compute
$\calP_{E_8}(x)$ to get
\[
 \calP_{E_8}(x) = x^8 + \delta_1x^7 + \delta_2x^6 + \delta_3x^5 + \delta_4x^4 + \delta_3x^3+\delta_2x^2+\delta_1x+ 1,
\]
where with $\ka = \ka_1 \ka_1\ka_2\ka_3\ka_4\ka_5\ka_6\ka_7\ka_8$,
\[
    \delta_1 = 8-\frac{g_1}{\kappa}, \quad \delta_2 = 28-6\frac{g_1}{\ka} + \frac{g_2}{\ka}, \quad \delta_3 = 
    56 -15\frac{g_1}{\ka}+ 4\frac{g_2}{\ka} -\frac{g_3}{\ka}, \quad \delta_4 = 70-20\frac{g_1}{\ka}+6\frac{g_2}{\ka}-2\frac{g_3}{\ka}+\frac{1}{\ka},
\]
\begin{align*}
    g_1 &= \ka_1  \ka_2  \ka_3  \ka_4  \ka_5  \ka_6 + \ka_1 \ka_2 \ka_3 \ka_4 \ka_5 \ka_8 + \ka_1 \ka_2 \ka_3 \ka_4 \ka_7 \ka_8 + 
 \ka_1 \ka_2 \ka_4 \ka_6 \ka_7 \ka_8 + \ka_1 \ka_2 \ka_5 \ka_6 \ka_7 \ka_8 \\
 & \;\;\;+ \ka_1 \ka_4 \ka_5 \ka_6 \ka_7 \ka_8 - 
 \ka_3 \ka_4 \ka_5  \ka_6  \ka_7  \ka_8, \\
 g_2 &= \ka_1 \ka_2 \ka_3 \ka_4 + \ka_1 \ka_2 \ka_4 \ka_6 + \ka_1 \ka_2 \ka_5 \ka_6 + \ka_1 \ka_4 \ka_5 \ka_6 + \ka_3 \ka_4 \ka_5 \ka_6 + \ka_1 \ka_2 \ka_4 \ka_8 + \ka_1 \ka_2 \ka_5 \ka_8 \\& \quad + \ka_1 \ka_4 \ka_5 \ka_8 + \ka_3 \ka_4 \ka_5 \ka_8 + \ka_1 \ka_2 \ka_7 \ka_8 + \ka_1 \ka_4 \ka_7 \ka_8 + \ka_3 \ka_4 \ka_7 \ka_8 
+ \ka_4 \ka_6 \ka_7 \ka_8  + \ka_5 \ka_6 \ka_7 \ka_8, \\
g_3 &= \ka_1 \ka_2 + \ka_1 \ka_4 + \ka_3 \ka_4 + \ka_4 \ka_6 + \ka_5 \ka_6  + \ka_4 \ka_8 + \ka_5 \ka_8 +\ka_7\ka_8.
\end{align*}
Using the explicit values of $(a_1, \ldots, a_8)$ in Table \ref{table:A}, we get 
$\delta_j = \gamma_j$ given in \leref{:D-poly} for $j \in [1, 4]$.

\bibliographystyle{alpha}
\bibliography{ref-July13-2025}
\end{document}